\documentclass{amsart}

\usepackage[all]{xy}
\usepackage{amssymb,amscd,epsfig}

\def\Z{{\mathbb Z}}
\def\Q{{\mathbb Q}}
\def\R{{\mathbb R}}
\def\C{{\mathbb C}}
\def\bH{{\mathbb H}}
\def\P{{\mathbb P}}
\def\cC{{\mathcal C}}
\def\cG{{\mathcal G}}
\def\F{{\mathcal F}}
\def\M{{\mathcal M}}
\def\O{{\mathcal O}}
\def\U{{\mathcal U}}

\def\d{\delta}
\def\e{\epsilon}
\def\w{\omega}
\def\g{\mathfrak g}
\def\D{\Delta}
\def\G{\Gamma}

\def\del{{\partial}}
\def\delbar{\overline{\partial}}

\def\wbar{\overline{w}}
\def\zbar{\overline{z}}
\def\sigmahat{\widehat{\sigma}}

\def\alphatilde{\tilde{\alpha}}
\def\B{\widehat{B}}
\def\E{E}
\def\W{\tilde{W}}
\def\Ebar{\overline{E}}

\def\Fdual{\check{\F}}
\def\Xbar{{\overline{X}}}

\def\dot{{\bullet}}
\def\bs{\backslash}
\def\blank{\phantom{x}}
\def\comp{~\widehat{\!}{\;}}

\def\sh{{\mathrm sh}}
\def\Ch{\mathit{Ch}}
\def\Ad{\mathit{Ad}}
\def\ev{\mathrm{ev}}
\def\H{{\mathrm{Hodge}}}

\def\HDR{H_{DR}}
\def\RTW{\R_{TW}^\dot}
\def\rtw{R_{TW}^\dot}
\def\simp{{\pmb{\Delta}}}
\def\Simp{{\mathrm{Simp}}}
\def\cosmp{[\dot]}
\def\Gm{{\mathbb G}_m}
\def\I{{I_\dot}}
\def\So{S^{(0)}}

\def\liminj#1{\lim_{\stackrel{\longrightarrow}{#1}}}
\def\limproj#1{\lim_{\stackrel{\longleftarrow}{#1}}}

\def\dlogz#1{\frac{d#1}{#1}}
\def\dlog1#1{\frac{d#1}{1-#1}}
\def\dlogq{\frac{d(xy)}{1-xy}}

\newcommand\id{\operatorname{id}}
\newcommand\im{\operatorname{im}}
\newcommand\Ln{\operatorname{ln}}
\newcommand\Hom{\operatorname{Hom}}
\newcommand\End{\operatorname{End}}
\newcommand\Ext{\operatorname{Ext}}
\newcommand\Jac{\operatorname{Jac}}
\newcommand\Pic{\operatorname{Pic}}

\newcommand\sgn{\operatorname{sgn}}
\renewcommand\Im{\operatorname{Im}}

\newcommand\Res{\operatorname{Res}}

\newtheorem{theorem}{Theorem}[section]
\newtheorem{lemma}[theorem]{Lemma}
\newtheorem{proposition}[theorem]{Proposition}
\newtheorem{corollary}[theorem]{Corollary}
\newtheorem{conjecture}[theorem]{Conjecture}

\theoremstyle{definition}
\newtheorem{example}[theorem]{Example}

\theoremstyle{remark}
\newtheorem{remark}[theorem]{Remark}

\begin{document}

\title[Iterated Integrals and Algebraic Cycles]
{Iterated Integrals and Algebraic Cycles: Examples and Prospects}

\author{Richard Hain}
\address{Department of Mathematics\\ Duke University\\
Durham, NC 27708-0320}
\email{hain@math.duke.edu}

\date{\today}

\thanks{Supported in part by grants from the National Science Foundation.}



\maketitle



The goal of this paper is to produce evidence for a connection between
the work of Kuo-Tsai Chen on iterated integrals and de~Rham homotopy theory
on the one hand, and the work of Wei-Liang Chow on algebraic cycles on the
other. Evidence for such a profound link has been emerging steadily since
the early 1980s when Carlson, Clemens and Morgan \cite{ccm} and Bruno Harris
\cite{harris} gave examples where the periods of non-abelian iterated
integrals coincide with the periods of homologically trivial algebraic
cycles. Algebraic cycles and the classical Chow groups are nowadays considered
in the broader arena of motives, algebraic $K$-theory and higher Chow groups.
This putative connection is best viewed in this larger context. Examples
relating iterated integrals and motives go back to Bloch's work
on the dilogarithm and regulators \cite{bloch} in the mid 1970s, which was
developed further by Beilinson \cite{beilinson} and Deligne (unpublished).
Further evidence to support a connection between de~Rham homotopy theory and
iterated integrals includes \cite{morgan,aomoto,hain:dht,wojtkowiak,pulte,
deligne:p1,hain-macp,zagier,yang,goncharov:trilog,hain-yang,shiho,terasoma,
colombo,zhao,cushman}. Chen would have been delighted by these developments,
as he believed iterated integrals and loopspaces contain non-trivial geometric
information and would one day become a useful mathematical tool outside
topology.

The paper is largely expository, beginning with an introduction to iterated
integrals and Chen's de~Rham theorems for loop spaces and fundamental groups.
It does contain some novelties, such as the de~Rham theorem for fundamental
groups of smooth algebraic curves in terms of ``meromorphic iterated integrals
of the second kind,'' and the treatment of the Hodge and weight filtrations
of the algebraic de~Rham cohomology of loop spaces of algebraic varieties in
characteristic zero. A generalization of the theorem of Carlson-Clemens-Morgan
in Section~\ref{cycles} is presented, although the proof is not complete for
lack of a rigorous theory of iterated integrals of currents. Even though
there is no rigorous theory, iterated integrals of currents are a useful
heuristic tool which illuminate the combinatorial and geometric content of
iterated integrals. The development of this theory should be extremely useful
for applications of de~Rham theory to the study of algebraic cycles. The
heuristic theory is discussed in Section~\ref{currents}.

A major limitation of iterated integrals and rational homotopy theory of
non-simply connected spaces is that they usually only give information about
nilpotent completions of topological invariants. This is particularly limiting
in many cases, such as when studying knots and moduli spaces of curves. By
using iterated integrals of twisted differential forms or certain convergent
infinite sums of iterated integrals, one may get beyond nilpotence.
Non-nilpotent iterated integrals and their Hodge theory should emerge when
studying the periods of extensions of variations of Hodge structure associated
to algebraic cycles in complex algebraic manifolds, when one spreads the
variety and the cycles. Some developments in the de~Rham theory, which
originate with a suggestion of Deligne, are surveyed in Section~\ref{beyond}.

Iterated integrals are the ``de~Rham realization'' of the cosimplicial
version of the cobar construction, a construction which goes back to Adams
\cite{adams}. The paper ends with an exposition of the cobar construction.
Logically, the paper could have begun with it, and some readers may prefer
to start there.  I hope that the examples in the paper will lead the reader
to the conclusion, first suggested by Wojtkowiak \cite{wojtkowiak}, that
the cosimplicial version of the cobar construction is important in algebraic
geometry, and that the numerous occurrences of iterated integrals as periods
of cycles and motives are not unrelated, but are the de~Rham manifestation of
a deeper connection between motives and the cobar construction.

This paper complements the survey article \cite{hain:bowdoin}, which
emphasizes the fundamental group. I highly recommend Chen's Bulletin article
\cite{chen:bams}; it surveys most of his work, and contains complete proofs
of many of his important theorems; it also contains a useful account of the
cobar construction. Polylogarithms are discussed from the point of view of
iterated integrals in \cite{hain:polylog}.

There is much beautiful mathematics that connects iterated integrals to
motives which is not covered in this paper. Most notable are Drinfeld's work
\cite{drinfeld}, in particular his associator, which appears in the study of
the motivic fundamental group of $\P^1 -\{0,1,\infty\}$, and the Kontsevich
integral \cite{kontsevich}, which appears in the construction of Vassiliev
invariants.
\bigskip

\noindent{\it Acknowledgements:} It is a great pleasure to acknowledge all
those who have inspired and contributed to my understanding of iterated
integrals, most notably Kuo-Tsai Chen, my thesis adviser, who introduced me to
them; Pierre Cartier, who influenced the way I think about them; and, Dennis
Sullivan who influenced me and many others through his seminal paper
\cite{sullivan} which still contains many paths yet to be explored.

\section{Differential Forms on Path Spaces}
\label{it_ints}

Denote the space of piece wise smooth paths $\gamma : [0,1] \to X$ in a 
smooth manifold $X$ by $PX$. Chen's iterated integrals can be defined
using any reasonable definition of differential form on $PX$, such at the
one used by Chen (see \cite{chen:bams}, for example). We shall denote the
de~Rham complex of $X$, $PX$, etc by $E^\dot(X)$, $E^\dot(PX)$, etc.

We will say that a function $\alpha : N \to PX$ from a smooth manifold
into $PX$ is {\it smooth} if the mapping
$$
\phi_\alpha : [0,1] \times N \to X
$$
defined by $(t,x) \mapsto \alpha(x)(t)$ is piecewise smooth in the sense
that there is a partition $0 = t_0 < t_1 < \dots < t_{n-1} < t_n = 1$ of
$[0,1]$ such that the restriction of $\phi_\alpha$ to each $[t_{j-1},t_j]
\times N$ is smooth.\footnote{Recall that a function $f:K \to \R$ from
a subset $K$ of $\R^N$ is smooth if there exists an open neighbourhood
$U$ of $K$ in $\R^N$ and a smooth function $g:U \to \R$ whose restriction
to $K$ is $f$.}

They key features of the de~Rham complex should satisfy are:
\begin{enumerate}
\item $E^\dot(PX)$ is a differential graded algebra;
\item if $N$ is a smooth manifold and $\alpha : N \to PX$ is smooth,
then there is an induced homomorphism
$$
\alpha^\ast : E^\dot(PX) \to E^\dot(N)
$$
of differential graded algebras;
\item if $D$ and $Q$ are manifolds and $D \times PX \to Q$ is smooth
(that is, $D\times N \to D\times PX \to Q$ is smooth for all smooth
$N \to PX$, where $N$ is a manifold), then there is an induced dga
homomorphism
$$
E^\dot(Q) \to E^\dot(D \times PX).
$$
\item If $D$ is compact oriented (possibly with boundary) of dimension $n$
and $p : D \times PX \to PX$ is the projection, then one has the integration
over the fiber mapping
$$
p_\ast : E^{k+n}(D\times PX) \to E^k(PX)
$$
which satisfies
$$
p_\ast d \pm d\,p_\ast = (p|_{\partial D})_\ast
$$
\end{enumerate}

Chen's approach is particularly elementary and direct. For him, a smooth
$k$-form on $PX$ is a collection
$w = (w_\alpha)$ of smooth $k$-forms, indexed by the smooth mappings
$\alpha : N_\alpha \to PX$, where $w_\alpha \in E^k(N_\alpha)$.
These are required to satisfy the following compatibility condition: if
$f : N_\alpha \to N_\beta$ is smooth, then
$$
w_\alpha = f^\ast w_\beta.
$$
Exterior derivatives are defined by setting $d(w_\alpha) = (dw_\alpha)$.
Exterior products are defined similarly. The de~Rham complex of $PX$
is a differential graded algebra.

This definition generalizes easily to other natural subspaces $W$ of $PX$,
such as loop spaces and fixed end point path spaces. Just replace $PX$ by
$W$ and consider only those $\alpha : N_\alpha \to PX$ that factor through
the inclusion $W \hookrightarrow PX$. It also generalizes to products of
such $W$ with a smooth manifold $Q$. To define a smooth form $w$ on
$Q\times W$, one need specify only the $w_\alpha$ for those smooth mappings
$\alpha$ of the form $\id \times \alpha: Q\times N \to Q\times W$.

Lest this seem ad hoc, I should mention that Chen developed an elementary
and efficient theory of ``differentiable spaces'', the category of which
contains the category of smooth manifolds and smooth maps, which is closed
under taking path spaces and subspaces. Each differentiable space has
a natural de~Rham complex which is functorial under smooth maps. The details
can be found in his Bulletin article \cite{chen:bams}.

\section{Iterated Integrals}

This is a brief sketch of iterated integrals. I have been deliberately
vague about the signs as they depend on choices of conventions which do not
play a crucial role in the theory. Another reason I have omitted them in
this discussion is that, by using different sign conventions from those of
Chen, I believe one should be able to make the signs in many formulas
conform more to standard homological conventions. Chen's sign conventions
are given in Theorem~\ref{alg_desc} and will be used in all computations
in this paper.

Suppose that $w_1,\dots,w_r$ are differential forms on $X$, all of positive
degree. The iterated integral
$$
\int w_1 w_2 \dots w_r
$$
is a differential form on $PX$ of degree
$-r + \deg w_1 + \deg w_2 + \cdots + \deg w_r$.
Up to a sign (which depends on one's conventions)
$$
\int w_1 w_2 \dots w_r =
\pi_\ast \phi^\ast
(p_1^\ast w_1 \wedge p_2^\ast w_2 \wedge \cdots \wedge p_r^\ast w_r)
$$
where
\begin{enumerate}
\item $p_j : X^r \to X$ is projection onto the $j$th factor,
\item
$\Delta^r = \{(t_1,\dots,t_r) : 0 \le t_1 \le t_2 \le \cdots \le t_r \le 1\}$
is the {\it time ordered} form of the standard $r$-simplex,
\item $\phi : \Delta^r \times PX \to X^r$
is the {\it sampling map}
$$
\phi(t_1,\dots,t_r,\gamma) = (\gamma(t_1),\gamma(t_2),\dots,\gamma(t_r)),
$$
\item $\pi_\ast$ denotes integration over the fiber of the projection
$$
\pi : \Delta^r \times PX \to PX.
$$
\end{enumerate}

When each $w_j$ is a 1-form, $\int w_1\dots w_r$ is a function $PX \to \R$.
Its value on the path $\gamma : [0,1] \to X$ is the {\it time ordered
integral}
\begin{equation}
\label{def}
\int_\gamma w_1 \dots w_r :=
\int_{0 \le t_1 \le t_2 \le \cdots \le t_n \le 1} f_1(t_1) \dots f_r(t_r)
dt_1 \dots d t_r,
\end{equation}
where $\gamma^\ast w_j = f_j(t) dt$. Iterated integrals of degree zero
are called {\it iterated line integrals}.

The space of iterated integrals on $PX$ is the subspace $\Ch^\dot(PX)$ of
its de~Rham complex spanned by all differential forms of the form
\begin{equation}
\label{term}
p_0^\ast w' \wedge p_1^\ast w'' \wedge \int w_1\dots w_r
\end{equation}
where for $a \in [0,1]$, $p_a : PX \to X$ is the evaluation at time $a$
mapping $\gamma \mapsto \gamma(a)$.

If $W$ is a subspace of $PX$ (such as a fixed end point path
space, the free loop space, a pointed loop space), we shall denote the
subspace of its de~Rham complex generated by the restrictions of iterated
integrals to it by $\Ch^\dot(W)$ and call it the {\it Chen complex} of $W$.
It is naturally filtered by length:
$$
\Ch_0^\dot(W) \subseteq \Ch_1^\dot(W) \subseteq \Ch_2^\dot(W)
\subseteq \cdots \subseteq \Ch^\dot(W),
$$
where $\Ch_s^\dot(W)$ consists of all iterated integrals that are sums
of terms (\ref{term}) where $r\le s$. 

The standard formula
$$
\pi_\ast d \pm d\,\pi_\ast = (\pi|_{\partial\Delta^r})_\ast
$$
implies that iterated integrals are closed under exterior differentiation
and that, with suitable signs (depending on one's conventions),
\begin{multline*}
d\int w_1\dots w_r = 
\sum_{j=1}^r \pm \int w_1\dots dw_j \dots w_r \cr
+ \sum_{j=1}^{r-1}
\pm \int w_1 \dots w_{j-1} (w_j\wedge w_{j+1}) w_{j+2} \dots w_r \cr
\pm  \bigg( \int w_1 \dots w_{r-1}\bigg) \wedge p_1^\ast w_r
\pm p_0^\ast w_1 \wedge \int w_2 \dots w_r.
\end{multline*}
This implies that each $\Ch_s^\dot(PX)$, and thus each $\Ch_s^\dot(W)$, is
closed under exterior differentiation.

The standard triangulation\footnote{This is 
$$
\Delta^r \times \Delta^s =
\bigcup_{\sigma \in \sh(r,s)}
\{(t_{\sigma(1)}, t_{\sigma(2)}, \dots, t_{\sigma(r+s)}): 0 \le t_1 \le
\dots \le t_r \le 1, 0 \le t_{r+1} \le \dots \le t_{r+s} \le 1\}
$$
}
of $\Delta^r \times \Delta^s$ gives the shuffle product formula
\begin{equation}
\int w_1 \dots w_r \wedge \int w_{r+1} \dots w_{r+s}
= \sum_{\sigma \in \sh(r,s)} \pm
\int w_{\sigma(1)} w_{\sigma(2)} \dots w_{\sigma(r+s)}
\end{equation}
\label{shuffle}
where $\sh(r,s)$ denotes the set of shuffles of type $(r,s)$ --- that is,
those permutations $\sigma$ of $\{1,2,\dots,r+s\}$ such that
\begin{multline*}
\sigma^{-1}(1) < \sigma^{-1}(2) < \cdots < \sigma^{-1}(r) \text{ and } \cr
\sigma^{-1}(s+1) < \sigma^{-1}(s+2) < \cdots < \sigma^{-1}(s+r).
\end{multline*}
With this product, $\Ch^\dot(W)$ is a differential graded algebra (dga).

In many applications, one considers the restrictions of iterated integrals
to the fixed end-point path spaces
$$
P_{x,y}X := \{\gamma \in PX : \gamma(0) = x, \gamma(1) = y\}.
$$
Multiplication of paths
$$
\mu : P_{x,y}X \times P_{y,z}X \to P_{x,z}X
$$
induces a map of the complex of iterated integrals:\footnote{Here we use
the convention that when $s=0$, $\int \phi_1 \dots \phi_s = 1$.}
$$
\mu^\ast \int w_1\dots w_r
= \sum_{j=1}^r
\pi_1^\ast \int w_1\dots w_j \wedge \pi_2^\ast \int w_{j+1} \dots w_r
$$
where $\pi_1$ and $\pi_2$ denote the projections onto the first and second
factors of $P_{x,y}X \times P_{y,z}X$. The inverse mapping
$$
P_{x,y}X \to P_{y,x}X,\qquad \gamma \mapsto \gamma^{-1}
$$
induces the ``antipode''
$$
\int w_1 \dots w_r \mapsto \pm \int w_r \dots w_1.
$$
The closed iterated line
integrals $H^0(\Ch^\dot(P_{x,y}X))$ are precisely those iterated line
integrals that are constant on homotopy classes of paths relative to 
their endpoints.

When $x=y$, the Chen complex of $P_{x,x}X$ is a differential
graded Hopf algebra with diagonal
$$
\int w_1 \dots w_r \mapsto \sum_{j=1}^r
\int w_1\dots w_j \otimes \int w_{j+1} \dots w_r.
$$
Its cohomology $H^\dot(Ch^\dot(P_{x,x}X))$ is a graded Hopf algebra with
antipode. Each element of $H^\dot(Ch^\dot(P_{x,x}X))$ defines a function
$\pi_1(X,x) \to \R$.

Restricting elements of $\Ch^\dot(P_{x,x}X)$ to the constant loop $c_x$ at
$x$ defines a natural augmentation
$$
Ch^\dot(P_{x,x}X) \to \R.
$$
Denote its kernel by $I\Ch^\dot(P_{x,x}X)$. These are the iterated integrals
on the loop space $P_{x,x}X$ ``with trivial constant term.''

Of course, if one takes iterated integrals of complex-valued forms
$E^\dot(X)_\C$, then one obtains complex-valued iterated integrals. We shall
denote the Chen complex of complex-valued iterated integrals by
$\Ch^\dot(PX)_\C$, $\Ch^\dot(P_{x,y}X)_\C$, etc.

\section{Loop Space de~Rham Theorems}
\label{derham}

Chen proved many useful de~Rham type theorems. (A comprehensive list can
be found in Section~2 of \cite{hain:dht}.) In this section we present
those of most immediate interest.

\begin{theorem}
\label{loop_dr}
If $X$ is a simply connected manifold, then integration induces a natural
Hopf algebra isomorphism
$$
H^\dot(\Ch^\dot(P_{x,x}X)) \cong H^\dot(P_{x,x}X;\R).
$$
\end{theorem}

This, combined with standard algebraic topology, gives a de~Rham theorem
for homotopy groups of simply connected manifolds. First a review of the
topology:

Suppose that $(Z,x)$ is a connected, pointed topological space and that $A$
is any coefficient ring. Consider the adjoint
$$
h^t : H^\dot(Z;A) \to \Hom_\Z(\pi_\dot(Z,z),A)
$$
of the Hurewicz homomorphism.
An element of $H^\dot(Z;A)$ is {\it decomposable} if it is in the image
of the cup product mapping
$$
H^{>0}(Z;A) \otimes H^{>0}(Z;A) \to H^\dot(Z;A).
$$
The set of indecomposable elements of the ring $H^\dot(Z;A)$ is defined
by
$$
QH^\dot(Z;A) := H^\dot(Z;A)/\text{\{the decomposable elements\}}.
$$
Since the cohomology ring of a sphere has no decomposables, the kernel of
$h^t$ contains the decomposable elements of $H^\dot(Z;A)$, and therefore
induces a mapping
$$
e : QH^\dot(Z;A) \to \Hom_\Z(\pi_\dot(Z,z),A).
$$
Typically, this mapping is far from being an isomorphism. However, if
$A$ is a field of characteristic zero and $Z$ is a connected $H$-space,
then $e$ is an isomorphism. (This is a Theorem of Cartan and Serre ---
cf.~\cite{milnor-moore}.)

When $X$ is simply connected, $P_{x,x}X$ is a connected $H$-space. Chen's
de~Rham theorem and the Cartan-Serre Theorem imply that
integration induces an isomorphism
$$
QH^j(P_{x,x}X;\R) \stackrel{\simeq}{\longrightarrow}
\Hom(\pi_j(P_{x,x},c_x),\R)  \cong  \Hom(\pi_{j+1}(X,x),\R)
$$
for each $j$.

There is a canonical {\it subcomplex} $Q\Ch^\dot(P_{x,x}X)$ of the Chen complex
of $P_{x,x}X$, which is isomorphic to the indecomposable iterated integrals
(see \cite{hain:indecomp}) and whose cohomology is
$QH^\dot(\Ch^\dot(P_{x,x}X))$. This and Chen's de~Rham Theorem above then
yield the following de~Rham Theorem for homotopy groups of simply connected
spaces:

\begin{theorem}
\label{derham_htpy}
Integration induces an isomorphism
$$
H^\dot(Q\Ch^\dot(P_{x,x}X)) \stackrel{\sim}{\longrightarrow}
\Hom(\pi_\dot(X,x),\R)
$$
of degree $+1$ of graded vector spaces.
\end{theorem}

Both sides of the display in this theorem are naturally ``Lie coalgebras.''
It is not difficult to show that the integration isomorphism respects
this structure.

We now turn our attention to non-simply connected spaces.
The augmentation $\epsilon : \Z\pi_1(X,x) \to \Z$ of the integral group
ring of $\pi_1(X,x)$ is defined by taking each element of the fundamental
group to 1. The augmentation ideal $J$ is the kernel of the augmentation
$\epsilon$. The diagonal mapping $\pi_1(X,x) \to \pi_1(X,x)\times \pi_1(X,x)$
induces a coproduct
$$
\Delta : \Z\pi_1(X,x)\to \Z\pi_1(X,x)\otimes \Z\pi_1(X,x).
$$

The most direct statement of Chen's de~Rham theorem for the fundamental
group is:

\begin{theorem}[Chen \cite{chen:tams}]
The integration pairing
$$
H^0(\Ch^\dot(P_{x,x}X)) \otimes \Z\pi_1(X,x) \to \C
$$
is a pairing of Hopf algebras under which
$H^0(\Ch_s^\dot(P_{x,x}X))$ annihilates $J^{s+1}$. The induced mapping
$$
H^0(\Ch_s^\dot(P_{x,x}X)) \to \Hom_\Z(\Z\pi_1(X,x)/J^{s+1},\C)
$$
is an isomorphism.
\end{theorem}

An elementary proof is given in \cite{hain:bowdoin}.
An equivalent statement, more amenable to generalization, will be given in
Section~\ref{beyond}. A second version uses the $J$-adic completion
$$
\R\pi_1(X,x)\comp := \limproj s \R\pi_1(X,x)/J^{s}
$$
of $\R\pi_1(X,x)$. It is a complete Hopf algebra (cf.\
\cite[Appendix~A]{quillen}), with diagonal
$$
\Delta : \R\pi_1(X,x)\comp \to
\R\pi_1(X,x)\comp\, \hat{\otimes} \R\pi_1(X,x)\comp
$$
induced by that of $\R\pi_1(X,x)$.

\begin{corollary}
If $H^1(X;\R)$ is finite dimensional, then integration induces a
natural homomorphism
$$
\R\pi_1(X,x)\comp \to \Hom(H^0(\Ch_s^\dot(P_{x,x}X)),\R)
$$
of complete Hopf algebras.
\end{corollary}

Of course, each of these theorems holds with complex coefficients if we
begin with complex-valued forms.

\begin{remark}
It is tempting to think that one can extend Chen's loop space de~Rham
theorem, or its homotopy version, to a de~Rham
theorem for higher homotopy groups of non-simply connected spaces. While
this is true for ``nilpotent spaces'' (which include Lie groups), it is most
definitely not true for most non-simply connected
spaces that one meets in day-to-day life. In fact, Example~\ref{limits}
shows that there is unlikely to be any reasonable
statement. One point we wish to make in this paper, however, is that in
arithmetic and algebraic geometry, the cohomology of iterated integrals is
intrinsic and may be a more interesting and geometric invariant of a complex
algebraic variety than its higher homotopy groups or loop space cohomology.
\end{remark}

\section{Multi-valued Functions}

In this section, we give several examples of interesting multi-valued
functions that can be obtained by integrating closed iterated line integrals.
Although elementary, these examples are reflections of the relationship between
iterated integrals and periods of certain canonical variations of mixed Hodge
structure.

The following result is easily proved by pulling back to the universal
covering of $X$ and using the definition (\ref{def}) of iterated line
integrals.

\begin{proposition}
\label{closed}
All closed iterated line integrals of length $\le 2$ on
$P_{x,y}X$ are of the form
$$
\sum_{j,k} a_{jk}\int \phi_j \phi_k + \int \xi + \text{ a constant}
$$
where each $\phi_j$ is a closed $1$-form on $X$, the $a_{jk}$ 
are scalars, and $\xi$ is a $1$-form on $X$ satisfying
$$
d\xi + \sum_{j,k} a_{jk}\, \phi_j \wedge \phi_k = 0.
$$
\end{proposition}

A {\it relatively closed} iterated integral is an element of $\Ch^\dot(PX)$
that is closed on $P_{x,y}X$ for all $x,y \in X$. The iterated line
integrals given by the previous result are relatively closed.

Multi-valued functions can be constructed by integrating relatively closed
iterated integrals. For example, suppose that $X$ is a Riemann surface and
that $w_1$ and $w_2$ are holomorphic differentials on $X$. Then
$$
\int w_1 w_2
$$
is closed on each $P_{x,y}X$. This means that for any fixed point $x_o \in X$,
the function
$$
x \mapsto \int_{x_o}^x w_1 w_2
$$
is a multi-valued function on $X$. It is easily seen to be holomorphic.

\begin{example}
If $X = \P^1(\C) - \{0,1,\infty\}$, then
$$
\int_0^x \frac{dz}{1-z}\frac{dz}{z}
$$
is a multi-valued holomorphic function on $X$. In fact, it is Euler's
dilogarithm, whose principal branch in the unit disk is defined by
$$
\Ln_2(x) = \sum_{n\ge 1}\frac{x^n}{n^2}.
$$
More generally, the $k$-logarithm
$$
\Ln_k(x) :=  \sum_{n\ge 1}\frac{x^n}{n^k}\qquad |x|< 1
$$
can be expressed as the length $k$ iterated integral
$$
\int_0^x \frac{dz}{1-z} \overbrace{\frac{dz}{z} \cdots \frac{dz}{z}}^{k-1}
$$
From this integral expression, it is clear that $\ln_k$ can be analytically
continued to a multi-valued function on $\C-\{0,1\}$.

Note that $\zeta(k)$, the value of the Riemann zeta function at an integer
$k > 1$, is $\ln_k(1)$. More information about iterated integrals and
polylogarithms can be found in \cite{hain:polylog}.
\end{example}

More generally, the multiple polylogarithms
$$
L_{m_1,\dots,m_n}(x_1,\dots,x_n) :=
\sum_{0<k_1<\cdots < k_n}
\frac{x_1^{k_1} x_2^{k_2} \dots x_n^{k_n}}{k_1^{m_1}k_2^{m_2}\dots k_n^{m_n}}
\qquad |x_j| < 1
$$
and their special values, Zagier's multiple zeta values $\zeta(n_1,\dots,n_m)$,
can be expressed as iterated integrals. For example,
$$
L_{1,1}(x,y) = \int_{(0,0)}^{(x,y)} \bigg(
\dlog1 y \dlog1 x + \dlogq\bigg(\dlog1 y - \dlog1 x - \dlogz x\bigg)
\bigg).
$$
This expression defines a well defined multi-valued function on
$$
\C^2 - \{(x,y) : xy(1-x)(1-y)(1-xy) \neq 0\}
$$
as the relation
$$
\dlog1 y \wedge \dlog1 x +
\dlogq \wedge \bigg(\dlog1 y - \dlog1 x - \dlogz x\bigg) = 0
$$
holds in the rational 2-forms on $\C^2$. Formulas for all 
multiple polylogarithms and other properties can be found in Zhao's paper
\cite{zhao}.

Closed iterated integrals that involve antiholomorphic 1-forms can also
yield multi-valued holomorphic functions.

\begin{proposition}
Suppose that $X$ is a complex manifold and $w_1,\dots, w_n$ are holomorphic
1-forms. If $\xi$ is a $(1,0)$ form such that
$$
\delbar \xi + \sum_{j,k} a_{jk} \wbar_j \wedge w_k = 0,
$$
then the multi-valued function
$$
F : x \mapsto
\sum_{j,k} a_{jk} \int_{x_o}^x \wbar_j w_k + \int_{x_o}^x \xi
$$
is well defined and holomorphic.
\end{proposition}

\begin{proof}
Since each $w_j$ is holomorphic,
$$
d F(x) = \sum_{j,k} a_{jk} \bigg(\int_{x_o}^x \wbar_j\bigg) w_k + \xi.
$$
Since $\xi$ has type $(1,0)$, $dF$ also has type $(1,0)$, which implies
that $F$ is holomorphic.
\end{proof}

\begin{example}
Take $X$ to be a punctured elliptic curve $E=\C/\Lambda-\{0\}$. The point
of this example is to show that the logarithm of the associated theta
function $\theta(z)$ is a twice iterated integral. We may
assume that $\Lambda = \Z + \Z\tau$ where $\tau$ has positive imaginary
part. Denote the homology classes of the images of the intervals $[0,1]$
and $[0,\tau]$ by $\alpha$ and $\beta$, respectively. These form a symplectic
basis of $H_1(E,\Z)$. The normalized abelian differential is $dz$ and
$dz = \alpha^\ast + \tau \beta^\ast$
from which it follows that
$$
d\zbar \wedge dz = 2i\Im \tau \in H^2(E,\C) \cong \C.
$$
The multi-valued differential $\mu = (z-\zbar)dz$ satisfies
$$
\mu(z+1) = \mu(z) \text{ and } \mu(z+\tau) = \mu(z) + 2i\Im\tau dz.
$$
On the other hand, the corresponding theta function $\theta(z) :=
\theta(z,\tau)$ satisfies
$$
\theta(z+1) = \theta(z) \text{ and } \theta(z+\tau) =
\exp(-i\pi\tau-2\pi iz)\theta(z).
$$
Thus
$$
\frac{d\theta}{\theta}(z+1) = \frac{d\theta}{\theta}(z) \text{ and }
\frac{d\theta}{\theta}(z+\tau) = \frac{d\theta}{\theta}(z) -2\pi i dz.
$$
It follows that
$$
\xi := \frac{\Im\tau}{\pi} \frac{d\theta}{\theta} + \mu(z)
$$
is a single-valued differential on $E-\{0\}$ of type $(1,0)$ having a
logarithmic singularity at 0 which satisfies
$$
d\zbar \wedge dz + d\xi = 0 \text{ in } E^2(E-\{0\}).
$$
(In fact, these properties characterize it up to a multiple of $dz$.)
It follows that $\int d\zbar dz + \xi$ is relatively closed, so that
$$
x \mapsto \int_{x_o}^x d\zbar\, dz + \xi
$$
is a  multi-valued holomorphic function on $E-\{0\}$. Applying the definition
of iterated integrals yields
$$
\log\theta(x) = \log\theta(x_o) +
\frac{\pi}{\Im \tau}
\bigg(
\int_{x_o}^x (d\zbar\, dz + \xi) - \frac{1}{2}\big(z(x)-\zbar(x_o)\big)^2
+ \frac{1}{2}\big(z(x_o)-\zbar(x_o)\big)^2
\bigg),
$$
where $z(x) = \int_0^x dz$.\footnote{In the language of Beilinson and
Levin \cite{beilinson-levin}, $\int dz$ is the elliptic logarithm for $E$,
and $\log \theta$ is the elliptic dilogarithm.}
This example generalizes easily to theta functions of several variables by
replacing $E$ by a principally polarized abelian variety $A$ and $E-\{0\}$
by $A - \Theta$, where $\Theta$ is its theta divisor.
\end{example}

\begin{remark}
This example can be developed further along the lines of Beilinson
\cite{beilinson} and Deligne's approach to the dilogarithm. (An exposition
of this, from the point of view of iterated integrals, can be found in
\cite{hain:polylog}.) Set
$$
G =
\begin{pmatrix}
1 & \C & \C \cr
0 & 1 & \C \cr
0 & 0 & 1
\end{pmatrix}
\text{ and }
F^0 G =
\begin{pmatrix}
1 & \C & 0 \cr
0 & 1 & 0 \cr
0 & 0 & 1
\end{pmatrix}.
$$
The Lie algebra $\g$ of $G$ is the Lie algebra of nilpotent upper triangular
matrices. Let
$$
w =
\begin{pmatrix}
0 & d\zbar & \xi \cr
0 & 0 & dz \cr
0 & 0 & 0
\end{pmatrix}
\in E^1(E-\{0\}) \otimes \g.
$$
This form is integrable: $dw + w \wedge w = 0$. It follows that the 
iterated integral
$$
T = \begin{pmatrix}
1 & \int d\zbar & \int d\zbar\, dz + \int \xi\cr
0 & 1 & \int dz \cr
0 & 0 & 1
\end{pmatrix}
$$
is relatively closed. (See \cite{chen:bams} or \cite{hain:bowdoin}.)
It therefore defines a homomorphism
$$
\theta : \pi_1(E-\{0\},x_o) \to G,
\quad \gamma \mapsto \langle T,\gamma\rangle
$$
which is the monodromy representation of the flat connection on
$(E-\{0\})\times G$ defined by $w$. There is thus a generalized Abel-Jacobi
mapping
$$
\nu : E-\{0\} \to \G \bs G / F^0 G
$$
that takes $x$ to $\langle T , \gamma \rangle$, where $\gamma$ is any
path in $E-\{0\}$ from $x_o$ to $x$. It is holomorphic as can be seen
directly using the formulas above.

Denote the center
$$
\begin{pmatrix}
1 & 0 & \C \cr
0 & 1 & 0 \cr
0 & 0 & 1
\end{pmatrix}
$$
of $G$ by $Z$. The quotient $\G\bs G /(F^0G\cdot Z)$ is naturally isomorphic
to $E$ itself and the corresponding projection
$$
\G \bs G / F^0 G \to E
$$
is a holomorphic $\C^\ast$-bundle. The formulas in the previous example
show that the section $\nu$ above is a non-zero multiple of the section
$\theta$ of the line bundle $L\to E$ associated to the divisor $0 \in E$:
$$
\xymatrix{
 & L - (0\mathrm{-section})\ar[d] \ar[r]^(0.55){\simeq} &
 {\G \bs G /F^0G} \ar[d] \cr
E-\{0\} \ar[r] \ar[ur]^{\theta} \ar[urr]^(0.3){\nu} &
E \ar[r]^(0.4){\simeq} & \G\bs G/(F^0G\cdot Z)
}
$$

This has an interpretation in terms of variations of MHS, which can be
extracted from \cite{hain-zucker}. The construction given here of the
second albanese mapping is special case of the direct construction given
in \cite{hain:albanese}. There is an analogous construction with a similar
interpretation where the pair $(E,0)$ is replaced by an abelian
variety and its theta divisor $(A,\Theta)$.
\end{remark}

\section{Harmonic Volume}
\label{harm_vol}

Bruno Harris \cite{harris} was the first to explicitly combine Hodge theory
and and (non-abelian) iterated integrals to obtain periods of algebraic cycles.
Suppose that $C$ is a compact
Riemann surface of genus 3 or more. Choose any base point $x_o \in C$.
Suppose that $L_1$, $L_2$, $L_3$ are three disjoint, non-separating simple
closed curves on $C$. Let $\phi_j$ be the harmonic representative
of the Poincar\'e dual of $L_j$, $j=1,2,3$. Since the curves are pairwise
disjoint, the product of any two of the $\phi_j$s vanishes in cohomology.
Thus there are 1-forms $\phi_{jk}$ such that
$$
d\phi_{jk} + \phi_j\wedge \phi_k = 0
$$
and $\phi_{jk}$ is orthogonal to the $d$-closed forms. These two conditions
characterize the $\phi_{jk}$. By Proposition~\ref{closed}, the iterated line
integral
$$
\int \phi_j \phi_k + \phi_{jk}
$$
is closed in $\Ch^\dot(P_{x_o,x_o}C)$.

Choose loops $\gamma_j$ ($j=1,2,3)$, based at $x_o$ and that are freely
homotopic to the $L_j$. Harris sets
$$
I(L_1,L_2,L_3) = \int_{\gamma_3} \phi_1 \phi_2 + \phi_{12}.
$$
He shows that this integral is independent of the choice of the base point
$x_o$ and that
$$
I(L_{\sigma(1)},L_{\sigma(2)},L_{\sigma(3)}) = \sgn(\sigma)I(L_1,L_2,L_3)
$$
for all permutations $\sigma$ of $\{1,2,3\}$.

One can also use the $\phi_j$ to imbed $C$ into the three torus
$T=\R^3/\Z^3$. Define $\Phi : C \to T$
by
$$
\Phi(x) = \int_{x_o}^x (\phi_1,\phi_2,\phi_3) \mod \Z^3.
$$
If the coordinates in $\R^3$ are $(z_1,z_2,z_3)$, then 
$\Phi^\ast dz_j = \phi_j$ for $j=1,2,3$. Since $H^2(T,\Z)$ is spanned by
the $dz_j\wedge dz_k$ and since
$$
\int_\Phi dz_j \wedge dz_k = \int_C \phi_j \wedge \phi_k = 0,
$$
the image of $C$ is homologous to zero in $T$. One can therefore find a
3-chain $\Gamma$ in $T$ such that $\del \Gamma = \Phi_\ast C$. Since $\Gamma$
is only well defined up to a 3-cycle, the volume of
$\Gamma$ is only well defined mod $\Z$. Harris's first main result is:

\begin{theorem}
The volume of $\Gamma$ is congruent to $I(L_1,L_2,L_3)$ mod $\Z$.
\end{theorem}

By an elementary computation, the span in $\Lambda^3 H_1(C,\Z)$ of classes
$L_1\wedge L_2 \wedge L_3$ is the kernel $K$ of the mapping
$$
\Lambda^3 H_1(C,\Z) \to H_1(C,\Z)
$$
defined by
$$
a\wedge b \wedge c \mapsto (a\cdot b)c + (b\cdot c) a + (c \cdot a)b.
$$
The harmonic volume $I$ thus determines a point in the compact torus
$\Hom(K,\R/\Z)$.

There is a lot more to this story --- it has a deep relationship to the
algebraic cycle $C_x - C_x^-$ in the jacobian of $C$. This is best 
explained in terms of the Hodge theory of the operator $\delbar$ rather
than $d$. This shall be sketched in Section~\ref{hodge}.

\section{Iterated Integrals of Currents}
\label{currents}

There is no rigorous theory of iterated integrals of currents, although such
a theory would be useful provided it is not too technical. The theory of
iterated integrals makes essential use of the algebra structure of the
de~Rham complex. The problem one encounters when trying to develop a theory
of iterated integrals of currents is that products of currents are only
defined when the currents being multiplied (intersected) are sufficiently
smooth (or sufficiently transverse). Nonetheless, this point of view is
useful, even it if is not rigorous. The paper \cite{hain:links} was an
attempt at making these ideas rigorous and using them to study links.

\begin{example}
In this example $X$ is the unit interval. Suppose
that $a_1, a_2,\dots, a_r$ are distinct points in the interior of the unit
interval. Set $w_j = \delta(t-a_j)dt$, where $\delta(t)$ denotes the Dirac
delta function supported at $t=0$. Let $\gamma : [0,1] \to X = [0,1]$ be
the identity path. Recall that $\Delta^r$ is the time ordered simplex
$$
\Delta^r = \{(t_1,t_2,\dots, t_r) : 0 \le t_1 \le \dots \le t_r \le 1\}.
$$
By definition,
\begin{align*}
\int_\gamma w_1 w_2 \dots w_r &=
\int_{\Delta^r}
\delta(t_1 - a_1)\delta(t_2-a_2)\dots \delta(t_r - a_r)\, dt_1 dt_2\dots dt_r 
\cr &=
\int_{\Delta^r} \delta_{(a_1,\dots,a_r)}(t_1,\dots,t_r)\, dt_1 dt_2\dots dt_r.
\end{align*}
Since the $a_j$ are distinct numbers satisfying $0<a_j<1$,
$$
\int_\gamma w_1 w_2 \dots w_r =
\begin{cases}
1 & \text{if } a_1 < a_2 < \cdots < a_r, \cr
0 & \text{otherwise.}
\end{cases}
$$
\end{example}

\begin{example}
More generally, suppose that $H_1, \dots, H_r$ are real hypersurfaces
in a manifold $X$, each with oriented (and thus trivial) normal bundle.
Suppose that $\gamma \in PX$ is transverse
to the union of the $H_j$ --- that is, the endpoints of $\gamma$ do not
lie in the union of the $H_j$ and $\gamma$ does not pass through any
singularity of their union. We can regard each $H_j$ as a current, which
we shall denote by $w_j$. For such a path $\gamma$ which is transverse to
$H_j$,
$$
\int_\gamma w_j = (H_j \cdot \gamma) := \text{ the intersection number of
$H_j$ with $\gamma$.}
$$
For simplicity, suppose that $\gamma$ passes
through each $H_j$ at most once, at time $t=a_j$, say. Then
$$
\gamma^\ast w_j = \epsilon_j \delta_j(t-a_j)
$$
where $\epsilon_j$ is 1 if $\gamma$ passes through $H_j$ positively at
time $a_j$, and $-1$ if it passes through negatively. By the previous
example,
$$
\int_\gamma w_1 w_2 \dots w_r =
\int_0^1 \gamma^\ast w_1 \dots \gamma^\ast w_r =
\begin{cases}
\epsilon_1 \epsilon_2 \dots \epsilon_r &
\text{if } a_1 < a_2 < \cdots < a_r, \cr
0 & \text{otherwise.}
\end{cases}
$$
\end{example}

\begin{center}
\begin{figure}[!ht]
\epsfig{file=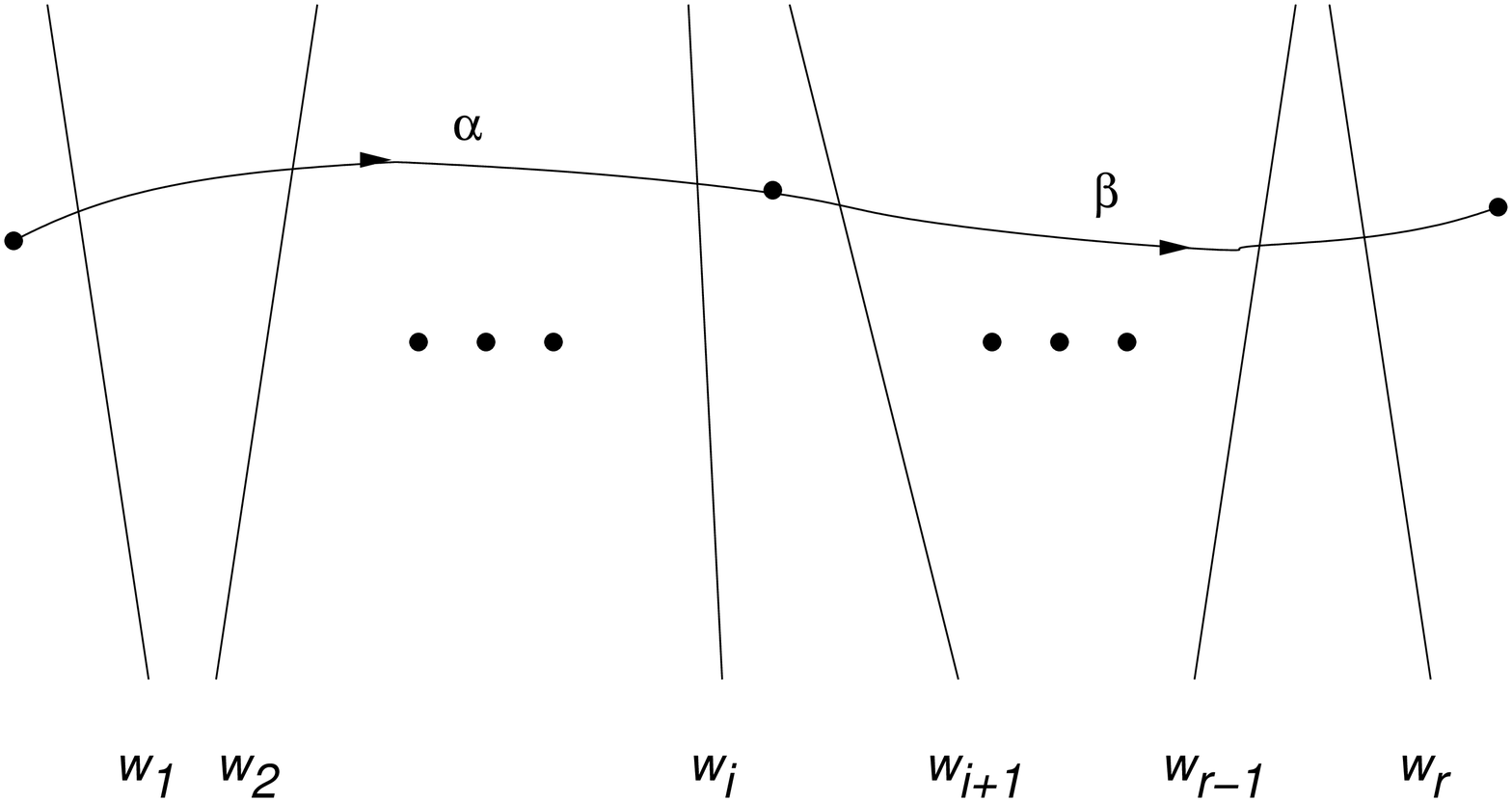, height=1.5in}
\caption{Pointwise product of iterated integrals}\label{prod}
\end{figure}
\end{center}

This formula can be used to give heuristic proofs of many basic properties
of iterated line integrals, such as the shuffle product formula, the
antipode, the coproduct and the differential. For example, suppose that
$w_1$, \dots, $w_r$ are 1-currents corresponding to oriented lines in the
plane and that $\alpha$ and $\beta$ are composable paths that are transverse
to the union of the supports of the $w_j$. (See Figure~\ref{prod}.)

Note that $\int w_1 \dots w_r$ is non-zero on $\alpha\beta$ if and only if
there is an $i$ such that $\alpha$ passes through $w_1, \dots, w_i$ in
order and $\beta$ passes through $w_{i+1}, \dots, w_r$ in order. In this
case
\begin{align*}
\int_{\alpha\beta} w_1\dots w_r &=
\int_\alpha w_i \dots w_i \int_\beta w_{i+1} \dots w_r \cr
&= \sum_{j=0}^r \int_\alpha w_j \dots w_i \int_\beta w_{j+1} \dots w_r
\end{align*}
as all the terms in the sum are zero except when $j=i$.

Examples using higher iterated integrals also exist. The simplest
I know of is a proof of the formula for the Hopf invariant
of a mapping $f:S^{4n-1} \to S^{2n}$. It is a nice exercise, using the
definition of iterated integrals, to show directly that if
$f : S^{4n-1} \to S^{2n}$ is smooth and $p$ and $q$ are distinct regular
values of $f$, then
$$
\big\langle \int \delta_p \delta_q, f\big\rangle
$$
is the linking number of $f^{-1}(p)$ and $f^{-1}(q)$ in $S^{4n-1}$.
Here $\delta_x$ denotes the $2n$-current associated to $x\in S^{2n}$.
This formula is equivalent to J.H.C.~Whitehead's integral formula for the
Hopf invariant \cite{whitehead} and Chen's version of it
\cite[p.~848]{chen:bams}.

\subsection{First steps}
\label{first_steps}
Suppose that $\G_1,\dots,\G_r$ are closed submanifolds of $X$ (possibly
with boundary), where $\G_j$ has codimension $d_j$. Denote
the $d_j$-current determined by $\G_j$ by $\d_j$. Suppose that $N$ is
a compact manifold and that $\alpha : N \to PX$ is smooth. We shall say that
$\alpha$ is {\it transverse to} $\int \d_1 \dots \d_r$ if the mapping
$$
\alphatilde : \D^r \times N \to X^r, \qquad
\big((t_1,\dots,t_r),u\big) \mapsto
\big(\alpha(u)(t_1),\dots,\alpha(u)(t_r)\big)
$$
is transverse to the submanifold $\G := \G_1 \times \cdots \times \G_r$ of
$X^r$. That is, the restriction of $\alphatilde$ to each stratum of
$\Delta^r \times N$ is transverse to each boundary stratum of $\G$.

If $N$ is has dimension $-r + d_1 + \dots + d_r$ and $\alphatilde$ is
transverse to $\G$, then we can evaluate the iterated integral
$$
\int \d_1 \dots \d_r
$$
on $\alpha$. This transversality condition is satisfied in each of the
examples above.

\section{The Reduced Bar Construction}
\label{bar}

Chen discovered that the iterated integrals on a smooth manifold have
a purely algebraic description \cite{chen:bams,chen:bar}. This algebraic
description is an important technical tool as it allows the computation
of various spectral sequences one obtains from iterated integrals,
applications to Hodge theory, and it facilitates the algebraic de~Rham
theory of iterated integrals for varieties over arbitrary algebraically
closed fields. (Cf.\ Section~\ref{alg_dr}.) This algebraic description is
expressed in terms of the {\it reduced bar construction}, a variant of the
more standard bar construction \cite{eilenberg-moore}, which is dual to
Adam's cobar construction \cite{adams}. Chen's version has the useful
property that it generates no elements of negative degree when applied
to a non-negatively graded dga with elements of degree zero, unlike the
standard version of the bar construction.

In this section, we use Chen's conventions for iterated integrals. In
particular, our description of the reduced bar construction gives a
precise formula for the exterior derivative of iterated integrals.

Suppose that $A^\dot$ is a differential graded algebra
(hereafter denoted dga) and that $M^\dot$ and $N^\dot$ are
complexes which are modules over $A^\dot$. That is, the
structure maps
$$
A^\dot \otimes M^\dot \to M^\dot \text{ and }
A^\dot\otimes N^\dot \to N^\dot
$$
are chain maps. We shall suppose that $A^\dot$, $M^\dot$ and $N^\dot$
are all non-negatively graded. Denote the subcomplex of $A^\dot$ consisting
of elements of positive degree by $A^{>0}$.

The {\it (reduced) bar construction} $B(M,A^\dot,N)$ is  defined as follows.
The underlying graded vector space is a quotient of the graded vector space
$$
T(M^\dot,A^\dot,N^\dot) :=
\bigoplus_s M^\dot \otimes\left(A^{>0}[1]^{\otimes r}\right) \otimes N^\dot.
$$
Following convention $m\otimes a_1\otimes \dots \otimes a_r \otimes n
\in T(M^\dot,A^\dot,N^\dot)$ will be denoted by $m[a_1|\dots|a_r]n$.
To obtain the vector space underlying the bar construction, we impose the
relations
\begin{align*}
m[dg|a_1|\dots|a_r]n & = m[ga_1|\dots|a_r]n - m\cdot g[a_1|\dots|a_r]n;\cr
\begin{split}
m[a_1|\dots|a_i|dg|a_{i+1}|\dots|a_r]n & = 
m[a_1|\dots|a_i|g\,a_{i+1}|\dots|a_r]n \cr
& \qquad - m[a_1|\dots|a_i\,g|a_{i+1}|\dots|a_r]n \quad 1\le i < s;
\end{split}\\
m[a_1|\dots|a_r|dg]n  & = m[a_1|\dots|a_r]g\cdot n - m[a_1|\dots|a_r\,g]n;\cr
m[dg]n & = 1 \otimes g\cdot n  - m\cdot g \otimes 1
\end{align*}
Here each $a_i \in A^{>0}$, $g\in A^0$, $m\in M^\dot$,  $n\in N^\dot$,
and $r$ is a positive integer.

Define an endomorphism $J$ of each graded vector space by 
$J: v\mapsto (-1)^{\deg v}v$. The differential is defined as 
$$
d = d_M\otimes 1_T \otimes 1_N + J_M\otimes d_B \otimes 1_N +
J_M \otimes J_T \otimes d_N + d_C.
$$
Here $T$ denotes the tensor algebra on $A^{>0}[1]$, $d_B$ is defined by
\begin{multline}
\label{diffl}
d_B[a_1|\dots|a_r] =
\sum_{1\le i \le r} (-1)^i [Ja_1|\dots|Ja_{i-1}|da_i|a_{i+1}|\dots|a_r]
\cr \hfill + \sum_{1 \le i < r}
(-1)^{i+1}[Ja_1|\dots|Ja_{i-1}|Ja_i\wedge a_{i+1}|a_{i+2}|\dots|a_r]
\end{multline}
and $d_C$ is defined by
$$
d_C m[a_1|\dots|a_r]n = (-1)^r
Jm[Ja_1|\dots|Ja_{r-1}]a_r \cdot n - Jm\cdot a_1 [a_2|\dots|a_r]n.
$$

The reduced bar construction $B(M^\dot,A^\dot,N^\dot)$ has a standard
filtration
$$
M^\dot\otimes N^\dot = B_0(M^\dot,A^\dot,N^\dot)
\subseteq B_1(M^\dot,A^\dot,N^\dot)
\subseteq B_2(M^\dot,A^\dot,N^\dot) \subseteq \cdots
$$
by subcomplexes, which is called the {\it bar filtration}. The subspace
$$
B_s(M^\dot,A^\dot,N^\dot)
$$
is defined to be the span of those $m[a_1|\dots|a_r]n$ with $r\le s$.
When $A^\dot$ has connected homology (i.e., $H^0(A^\dot) = \R$), the
corresponding (second quadrant) spectral sequence, which is called the
{\it Eilenberg-Moore spectral sequence} (EMss), has $E_1$ term
$$
E_1^{-s,t} =
\left[M^\dot\otimes H^{>0}(A^\dot)^{\otimes s}\otimes N^\dot\right]^t.
$$
A proof can be found in \cite{chen:bar}. This computation has the following
useful consequence:

\begin{lemma}
\label{qism}
Suppose that $A_j^\dot$ is a dga and $M_j^\dot$ and $N_j^\dot$ are right
and left $A_j^\dot$-modules, where $j=1,2$. Suppose that
$f_A :A_1\dot \to A_2^\dot$ is a dga homomorphism and
$$
f_M:M_1^\dot \to M_2^\dot \text{ and }f_N: N_1^\dot \to N_2^\dot
$$
are chain maps compatible with the the actions of $A_1^\dot$ and $A_2^\dot$
and $f$. If $f_A$, $f_M$ and $f_N$ induce isomorphisms on homology, then
so do the induced mappings
$$
B_s(M_1^\dot,A_1^\dot,N_1^\dot) \to B_s(M_2^\dot,A_2^\dot,N_2^\dot)
\text{ and }
B(M_1^\dot,A_1^\dot,N_1^\dot) \to B(M_2^\dot,A_2^\dot,N_2^\dot)
$$
\end{lemma}

When $A^\dot$, $M^\dot$ and $N^\dot$ are commutative dgas (in the graded
sense), and when the
$A^\dot$-module structure on $M^\dot$ and $N^\dot$ is determined by dga
homomorphism $A^\dot \to M^\dot$ and $A^\dot \to N^\dot$,
$B(M^\dot,A^\dot,N^\dot)$ is also a commutative dga. The multiplication is
given by the shuffle product:
\begin{multline*}
\big(m'[a_1|\dots|a_r]n'\big) \wedge \big(m''[a_{r+1}| \dots | a_{r+s}]n''\big)
\cr
= \sum_{\sh(r,s)} \pm 
(m'\wedge m'')[a_{\sigma(1)}| \dots | a_{\sigma(r+s)}](n'\wedge n'').
\end{multline*}
It is important to note that the shuffle product does not commute with the
differential when $A^\dot$ is not commutative.

Many complexes of iterated integrals may be described in terms of
reduced bar constructions of suitable triples. Here we give just one
example --- the iterated integrals on $P_{x,y}X$. A more complete list
of such descriptions can be found in \cite{chen:bams} and \cite[\S2]{hain:dht}.

Suppose that $X$ is a manifold and that $x_0$ and $x_1$ are points
of $X$. Evaluating at $x_j$, we obtain an augmentation
$\e_j : E^\dot(X)\to \R$ for $j=0,1$. Suppose that $A^\dot$ is a sub dga
of $E^\dot(M)$ and that both augmentations restrict to non-trivial
homomorphisms $\e_j : A^\dot \to \R$. We can take $M^\dot$ and $N^\dot$ both
to be $\R$, where the action is given by $\e_0$ and $\e_1$, respectively. Now
form the corresponding bar construction $B(\R,A^\dot,\R)$.

Define $\Ch^\dot(P_{x_0,x_1}X;A^\dot)$ to be the subcomplex of 
$\Ch^\dot(P_{x_0,x_1}X)$ spanned by those iterated integrals
$\int w_1\dots w_r$ where each $w_j \in A^\dot$.

\begin{theorem}
\label{alg_desc}
Suppose that $X$ is connected. If $H^0(A^\dot) \cong \R$ and the natural
map $H^1(A^\dot) \to H^1(X;\R)$ induced by the inclusion of $A^\dot$ into
$E^\dot(X)$ is injective, then the natural mapping
$$
B(\R,A^\dot,\R) \to \Ch^\dot(P_{x_0,x_1}X;A^\dot),
\qquad [w_1| \dots | w_r] \mapsto \int w_1 \dots w_r
$$
is a well defined isomorphism of differential graded algebras.
\end{theorem}

This and Adams' work \cite{adams} are the basic ingredients in the proof of
the loop space de~Rham theorem, Theorem~\ref{loop_dr}. The previous result
has many useful consequences, such as:

\begin{corollary}
\label{A_ints}
If $X$ is connected and $A^\dot$ is a sub dga of $E^\dot(X)$ for which
the inclusion $A^\dot \hookrightarrow E^\dot(X)$ induces an isomorphism
on homology, then the inclusion
$$
\Ch^\dot(P_{x,y}X;A^\dot) \hookrightarrow \Ch^\dot(P_{x,y}X)
$$
induces an isomorphism on homology.
\end{corollary}

This is proved using the previous two results. It has many uses, such
as in the next example, where it simplifies computations, and in Hodge
theory, where one takes $A^\dot$ to be the subcomplex of $C^\infty$
logarithmic forms when $X$ is the complement of a normal crossings divisor
in a complex projective algebraic manifold.

\begin{example}
A nice application of the results so far is to compute the loop space
cohomology $H^\dot(P_{x,x}S^n;\R)$ and real homotopy groups
$\pi_n(S^n,x)\otimes \R$ of the $n$-sphere ($n\ge 2$). This computation
is classical.

The first thing to do is to replace the de~Rham complex of $S^n$ by a 
sub dga $A^\dot$ which is as small as possible, but which computes the
cohomology of the sphere. To do this, choose an $n$-form $w$ 
whose integral over $S^n$ is 1 and take $A^\dot$ to be the dga consisting
of the constant functions and the constant multiplies of $w$. By
Corollary~\ref{A_ints}, the iterated integrals constructed from elements of
$A^\dot$ compute the cohomology of $S^n$. But these are all linear
combinations of
$$
\theta_m := \int \overbrace{w \dots w}^m, \quad m \ge 0.
$$
Each of these is closed, and no linear combination of them is exact. It
follows that
$$
H^j(P_{x,x}S^n;\R) \cong
\begin{cases}
\R\,\theta_m & j = m(n-1);\cr
0 & \text{ otherwise.}
\end{cases}
$$
The ring structure is also easily determined using the shuffle product 
formula (\ref{shuffle}). When $n$ is odd, we have $\theta_1^m = m! \theta_m$;
when $n$ is even
$$
\theta_1 \wedge \theta_{2m} = \theta_{2m+1}, \quad
\theta_1\wedge \theta_{2m+1} = 0, \text{ and }
\theta_2 \wedge \theta_{2m} = (m+1)\theta_{2m+2}.
$$

Applying Theorem~\ref{derham_htpy} we have:
$$
\pi_j(S^n,x)\otimes \R =
\begin{cases}
\R & j = n; \cr
\R & j = 2n-1 \text{ and $n$ even}; \cr
0 &  \text{ otherwise.}
\end{cases}
$$\qed
\end{example}

\begin{example}
\label{limits}
This example illustrates the limits of the ability of iterated integrals
to compute homotopy groups.\footnote{Minimal models do no better or
worse. If $(X,x)$ is a pointed topological space with minimal model
$\M_X^\dot$, there is a canonical Lie coalgebra isomorphism
$Q\M_X^\dot \cong H^\dot(\Ch^\dot(P_{x,x}X))$. This follows from
\cite[\S3]{chen:circular}.}
The main point is that there are continuous maps
$f:X \to Y$ between spaces that induce an isomorphism on cohomology, but
not on homotopy. Properties of the bar construction (cf.\ Lemma~\ref{qism})
imply that for such $f$ the mapping
$$
f^\ast : H^\dot(\Ch^\dot(P_{f(x),f(x)}Y)) \to H^\dot(\Ch^\dot(P_{x,x}X))
$$
is an isomorphism.

The prototype of such continuous functions is the mapping $X \to X^+$ from
a connected topological space $X$, with the property that the commutator
subgroup of $\pi_1(X,x)$ is perfect, to $X^+$, Quillen's plus construction.

By a standard trick, one can extend de~Rham theory (and hence iterated
integrals) to arbitrary topological spaces.\footnote{Basically, one replaces
a space by the simplicial set consisting of its singular chains. This is
canonically weak homotopy equivalent to the original space. One then can
work with the Thom-Whitney de~Rham complex of this simplicial set. It
computes the cohomology of the space and is functorial under continuous
maps.}
In this setting, one can take a perfect group $\G$ and consider the mapping
$$
\phi : B\G \to B\G^+
$$
from the classifying space of $\G$ to its plus construction. This mapping
induces an isomorphism on homology, and therefore a quasi-isomorphism
$$
\phi^\ast : E^\dot(B\G^+) \to E^\dot(B\G).
$$
This induces, by Lemma~\ref{qism}, an isomorphism
$$
H^\dot(\Ch^\dot(P_{x,x}B\G^+)) \to H^\dot(\Ch^\dot(P_{x,x}B\G))
$$
Since the universal covering of $B\G$ is contractible, $P_{x,x}B\G$
is a disjoint union of contractible sets indexed by the elements of $\G$.
On the other hand, $B\G^+$ is a simply connected $H$-space, the loop space
de~Rham theorem holds for it. It follows that
$$
QH^j(\Ch^\dot(P_{x,x}B\G)) \cong \Hom(\pi_{j+1}(B\G^+,x), \R).
$$

In particular, take $\G = SL(\Z)$, a perfect group. From Borel's work
\cite{borel}, we know that
$$
\pi_j(BSL(\Z)^+,x)\otimes \R \cong
\begin{cases}
\R & j \equiv 3 \mod 4 \cr
0 & \text{ otherwise.}
\end{cases}
$$

For those who would prefer an example with manifolds, one can approximate
$BSL(\Z)$ by a finite skeleton of $BSL_n(\Z)$ for some $n\ge 3$ or take
$\G$ to be a mapping class group in genus $g\ge 3$.
\end{example}

\subsection{An integral version}
\label{integral}
Suppose that $X$ is a topological space and that $R$ is a ring.
Each point $x$ of $X$ induces an augmentation $\e_x : S^\dot(X;R) \to R$
on the $R$-valued singular chain complex of $X$. If $x,y \in X$, we
have augmentations
$$
\e_x : S^\dot(X;R) \to R \text{ and } \e_y : S^\dot(X;R) \to R,
$$ 
which give $R$ two structures as a module over the singular cochains. We
can thus form the reduced bar construction $B(R,S^\dot(X;R),R)$.

The following result, which will be further elaborated in Section~\ref{cobar}
and is proved using Adams cobar construction, is needed to put an integral
structure on the cohomology of $\Ch_s^\dot(P_{x,y}X)$, regardless of whether
$X$ is simply connected or not.

\begin{proposition}[Chen \cite{chen:bams}]
\label{integrality}
For all $s \ge 0$, there are canonical isomorphisms
$$
H^\dot(B_s(\Z,S^\dot(X,\Z),\Z))\otimes_\Z \R \cong
H^\dot(B_s(\R,S^\dot(X,\R),\R)) \cong H^\dot(\Ch_s^\dot(P_{x,y}X)).
$$
\end{proposition}

It is very important to note that the na\"{\i}ve mapping
$$
B(I):B(\R,E^\dot(X),\R) \to B(\R,S^\dot(X;\R),\R), \quad
[w_1|\dots | w_r] \mapsto [I(w_1)| \dots | I(w_r)]
$$
induced by the integration mapping $I:E^\dot(X) \to S^\dot(X;\R)$, is
{\it not} a chain mapping. This is because $I$ is not an algebra
homomorphism (except in trivial cases), which implies that $B(I)$ is
not, in general, a chain mapping.

\section{Exact Sequences}

The algebraic description of iterated integrals gives rise to several
exact sequences useful in topology and Hodge theory. We shall concentrate
on iterated integrals of length $\le 2$ as this is the first interesting
case --- $H^k(I\Ch_1^\dot(P_{x,x}X))$ is just $H^{k+1}(X;\R)$.

\begin{lemma}
\label{ex_seq}
If $X$ is a connected manifold, then the sequence
\begin{multline*}
0 \to QH^{2d-1}(X;\R) \to H^{2d-2}(I\Ch_2^\dot(P_{x,x}X)) \to
[H^{>0}(X;\R)^{\otimes 2}]^{2d}\cr
\stackrel{\text{cup}}{\longrightarrow} H^{2d}(X;\R) \to QH^{2d}(X;\R) \to 0
$$
\end{multline*}
is exact.
This sequence has a natural $\Z$-form and exactness holds over $\Z$ as
well.
\end{lemma}

\begin{proof}[Sketch of Proof]
By the algebraic description of iterated integrals given in the previous
section, the sequence
$$
0 \to I\Ch_1^\dot(P_{x,x}X) \to I\Ch_2^\dot(P_{x,x}X) \to
(E^{>0}(X)/dE^0(X))^{\otimes 2} \to 0
$$
is exact. This gives rise to a long exact sequence. The formula for the
differential and the identification of $I\Ch_1^\dot(P_{x,x}X)$ with
$E^{>0}(X)/dE^0(X)$ imply that the connecting homomorphism is the cup product
$$
[H^{>0}(X;\R)^{\otimes 2}]^k \to H^k(X;\R). 
$$
The integrality statement follows from Prop.~\ref{integrality} using the
integral version of the reduced bar construction.
\end{proof}

Combining it with the de~Rham Theorems yields the following two results.
For the first, note that the function
$$
\pi_1(X,x) \to J/J^2, \qquad \gamma \mapsto (\gamma -1) + J^2
$$
is a homomorphism and induces an isomorphism
$$
H_1(X,\Z) \cong J/J^2.
$$
Here $J$ denotes the augmentation ideal of $\Z\pi_1(X,x)$.

\begin{corollary}
\label{ex_pi1}
For all connected manifolds $X$, the sequence
$$
0 \to H^1(X;\Z) \to \Hom(J/J^3,\Z) \stackrel{\psi}{\to} H^1(X;\Z)^{\otimes 2}
\stackrel{\text{cup}}{\longrightarrow} H^2(X;\Z)
$$
is exact. The mapping $\psi$ is dual to the multiplication mapping
$$
H_1(X;\Z)^{\otimes 2} \cong (J/J^2)^{\otimes 2} \to J/J^3.
$$
\end{corollary}

The analogue of this in the simply connected case is:

\begin{corollary}
\label{ex_pi3}
If $X$ is simply connected, then the sequences
$$
0 \to H^3(X;\Q) \to \Hom(\pi_3(X,x),\Q) \to S^2 H^2(X;\Q)
\stackrel{\text{cup}}{\longrightarrow} H^4(X)
$$
and
$$
0 \to H^3(X;\Z) \to H^2(P_{x,x}X;\Z) \to H^2(X;\Z)^{\otimes 2}
\stackrel{\text{cup}}{\longrightarrow} H^4(X;\Z)
$$
are exact.
\end{corollary}

\section{Hodge Theory}
\label{hodge}

Just as in the case of ordinary cohomology, Chen's de~Rham theory is much
more powerful when combined with Hodge theory, and is especially fertile
when applied to problems in algebraic geometry. The Hodge theory of
iterated integrals is best formalized in terms of Deligne's mixed Hodge
theory. I will not review Deligne's theory here, but (at the peril of 
satisfying nobody) will attempt to present the ideas in a way that will
make sense both the novice and the expert. More details can be found in
\cite{hain:dht,hain:hodge,hain:bowdoin,hain-zucker}.

\subsection{The riemannian case}
In the classical case, the Hodge theorem asserts that every
de~Rham cohomology class on a compact riemannian manifold has a unique
harmonic representative which depends, in general, on the metric. If
$X$ is a compact riemannian manifold, then every element of
$H^\dot(\Ch_s^\dot(P_{x,y}X))$ has a natural representative, which I shall
call ``harmonic'' even though  I do not know if it is annihilated by any
kind of laplacian on $P_{x,y}X$.

This is illustrated in the case $s=2$. Every closed iterated integral of
length $\le 2$ is of the form
\begin{equation}
\label{rep}
\sum_{j,k} a_{jk} \int w_j w_k + \xi
\end{equation}
where
$$
d \xi = \sum_{j,k} (-1)^{\deg w_j} a_{jk}\, w_j \wedge w_k.
$$
The iterated integral (\ref{rep}) is defined to be harmonic if each $w_j$ is 
harmonic and $\xi$ is co-closed (i.e., orthogonal to the closed forms). This
definition generalizes to iterated integrals of arbitrary length.
Classical harmonic theory on $X$ and the Eilenberg-Moore spectral sequence
imply that every element of $H^\dot(\Ch^\dot(P_{x,y}X))$ has a unique
harmonic representative.

Harris's work on harmonic volume (Section~\ref{harm_vol}) is a particularly
nice application of harmonic iterated integrals.

\subsection{The K\"ahler case}
This na\"{\i}ve picture generalizes to the case when $X$ is compact
K\"ahler. In classical Hodge theory, certain aspects of the Hodge Theorem,
such as the Hodge decomposition of the cohomology, are independent of the
metric. Similar statements hold for iterated integrals: specifically,
$H^\dot(\Ch^\dot(P_{x,y}X))$ has a natural mixed Hodge structure, the
key ingredient of which is the {\it Hodge filtration}, whose definition
we now recall.

The Hodge filtration
$$
E^\dot(X)_\C = F^0 E^\dot(X) \supseteq F^1 E^\dot(X) \supseteq E^\dot(X) \supseteq \cdots
$$
of the de~Rham complex of a complex manifold is defined by
$$
F^p E^\dot(X) := \bigoplus_{r \ge p} E^{r,\dot}(X)
= \left\{\parbox{2.5in}{differential forms for which each term of each
local expression has at least $p$ $dz$'s}\right\}
$$
where $E^{p,q}(X)$ denotes the differential forms of type $(p,q)$ on $X$.
Each $F^pE^\dot(X)$ is closed under exterior differentiation.

A fundamental consequence of the Hodge theorem for compact K\"ahler manifolds
is the following:

\begin{proposition}
\label{strictness}
If $X$ is a compact K\"ahler manifold, then the mapping
$$
H^\dot(F^pE^\dot(X)) \to H^\dot(X;\C)
$$
is injective and has image
$$
F^pH^\dot(X) := \bigoplus_{s\ge p} H^{s,t}(X)
$$
In other words, every class in $F^pH^\dot(X)$ is represented by a class
in $F^pE^\dot(X)$, and if $w \in F^p E^\dot(X)$ is exact in $E^\dot(X)_\C$,
one can find $\psi \in F^p E^\dot(X)$ such that $d\psi = w$.
\end{proposition}

The Hodge filtration extends naturally to complex-valued iterated integrals:
$F^p\Ch_s^\dot(P_{x,y}X)$ is the span of
$$
\int w_1 \dots w_r
$$
where $r\le s$ and $w_j \in F^{p_j}E^\dot(X)$, where
$p_1 + \dots + p_r \ge p$.
The weight filtration is simply the filtration by length:
$$
W_m\Ch^\dot(P_{x,y}X) = Ch_m^\dot(P_{x,y}X).
$$

The Hodge theory of iterated integrals for compact K\"ahler manifolds is
summarized in the following result. A sketch of a proof can be found in
\cite{hain:hodge} and a complete proof in \cite{hain:dht}.

\begin{theorem}
If $X$ is a compact K\"ahler manifold, then $\Ch^\dot(P_{x,y}X)$, endowed
with the Hodge and weight filtrations above, is a mixed Hodge complex. In
particular:
\begin{enumerate}
\item $H^\dot(F^p\Ch^\dot(P_{x,y}X)) \to H^\dot(\Ch^\dot(P_{x,y}X)_\C)$ is
injective;
\item $H^\dot(\Ch^\dot(P_{x,y}X))$ has a natural mixed Hodge structure with
Hodge and weight filtrations defined by
$$
F^p H^\dot(\Ch^\dot(P_{x,y}X)) = H^\dot(F^p\Ch^\dot(P_{x,y}X))
$$
and
$$
W_m H^k(\Ch^\dot(P_{x,y}X)) =
\im \big\{H^k(\Ch_{m-k}^\dot(P_{x,y}X) \to H^k(\Ch^\dot(P_{x,y}X))\big\}
$$
\end{enumerate}
If $H^1(X;\Q)=0$, this mixed Hodge structure is independent of the base point
$x$
\end{theorem}

This theorem generalizes to all complex algebraic manifolds (using logarithmic
forms) and to singular complex algebraic varieties (using simplicial methods).
Details can be found in \cite{hain:dht}.

\begin{corollary}
\label{mhs_exact}
If $X$ is a complex algebraic manifold, then $H^{2d-2}(I\Ch_2^\dot(P_{x,x}X))$
has a canonical mixed Hodge structure defined over $\Z$ and the sequence
\begin{multline*}
0 \to QH^{2d-1}(X) \to H^{2d-2}(I\Ch_2^\dot(P_{x,x}X)) \to
[H^{>0}(X)^{\otimes 2}]^{2d}\cr
\stackrel{\text{cup}}{\longrightarrow} H^{2d}(X) \to QH^{2d}(X) \to 0.
\end{multline*}
is exact in the category of $\Z$ mixed Hodge structures.
\end{corollary}

The minimal model approach to the Hodge theory for complex algebraic manifolds
was developed by Morgan in \cite{morgan}. From the point of view of Hodge
theory, iterated integrals have the advantage that they provide a rigid
invariant on which to do Hodge theory, whereas
the minimal model of a manifold is unique only up to a homotopy class of
isomorphisms, which makes the task of putting a mixed Hodge structure on a
minimal model more difficult. Chen's theory is also better suited to studying
the non-trivial role of the base point $x\in X$ in the theory, which is
particularly important when studying the Hodge theory of the fundamental group.
On the other hand, minimal models (and other non-rigid models) are an essential
tool in understanding how Hodge theory restricts fundamental groups and
homotopy types of complex algebraic varieties, as is illustrated by Morgan's
remarkable examples in \cite{morgan}.

\section{Applications to Algebraic Cycles}
\label{cycles}

Recall that a Hodge structure $H$ of weight $m$ consists of a finitely
generated abelian group $H_\Z$ and a bigrading
$$
H_\C = \bigoplus_{p+q = m} H^{p,q}
$$
of $H_\C = H_\Z\otimes \C$ by complex subspaces satisfying
$
H^{p,q} = \overline{H^{q,p}}\ .
$
The standard example of a Hodge structure of weight $m$ is the $m$th
cohomology of a compact K\"ahler manifold. Its dual, $H_m(X)$, is a Hodge
structure of weight $-m$.\footnote{Just define $H_m(X)^{-p,-q}$ to be the
dual of $H^{p,q}(X)$.}

For and integer $d$, the {\it Tate twist} $H(d)$ of $H$ is defined to be
the Hodge structure with the same underlying lattice $H_\Z$ but whose
bigrading has been reindexed:
$$
H(d)^{p,q} = H^{p+d,q+d}.
$$
Equivalently, $H(d)$ is the tensor product of $H$ with the 1-dimensional
Hodge structure $\Z(d)$ of weight $-2d$.

The category of Hodge structures is abelian, and closed under tensor products
and taking duals.

\subsection{Intermediate jacobians and Griffiths' construction}
The $d$th {\it intermediate jacobian} of a compact K\"ahler manifold
$X$ is defined by
$$
J_d(X) := J(H_{2d+1}(X)(-d)) \cong \Hom(F^{d+1}H^{2d+1}(X),\C)/H_{2d+1}(X;\Z).
$$
It is a compact, complex torus. For example, $J_0(X)$ is the albanese of $X$
and $J_{\dim X-1}(X)$ is $\Pic^0 X$, the group of isomorphism classes of
topologically trivial holomorphic line bundles over $X$.

Suppose $Z$ is an algebraic $d$-cycle in $X$, that is trivial in homology.
We can write $Z$ as the boundary of a $(2d+1)$-chain $\G$, which determines
a point $\int_\G$ in
$$
\Hom(F^{d+1}H^{2d+1}(X),\C)/H_{2d+1}(X;\Z) \cong J_d(X)
$$
by integration:
$$
\int_\G : [w] \mapsto \int_\G w
$$
where $w \in F^{d+1}E^{2d+1}(X)$. This mapping is well defined by Stokes'
Theorem, Proposition~\ref{strictness}, and because $F^{d+1}E^{2d}(Z) = 0$.

It is also convenient to define $J^d(X) = J_{n-d}(X)$, where $n$ is the
complex dimension of $X$.

\subsection{Extensions of mixed Hodge structures}
In this paragraph, we review some elementary facts about extensions of
mixed Hodge structures (MHS). Complete details can be found in \cite{carlson}.
Suppose that $A$ and $B$ are Hodge structures of weights $n$ and $m$,
respectively, and that
\begin{equation}
\label{ext}
0 \to B \to E \stackrel{\pi}{\to} A \to 0
\end{equation}
is an exact sequence of mixed Hodge structures. In concrete terms, this
means:
\begin{enumerate}
\item there is an exact sequence
\begin{equation}
\label{extn}
0 \to B_\Z \to E_\Z \stackrel{\pi}{\to} A_\Z \to 0;
\end{equation}
of finitely generated abelian groups;
\item $E_\C := E_\Z\otimes \C$ has a filtration
$\cdots \supseteq F^p E \supseteq F^{p+1} E \supseteq \cdots$
satisfying
$$
B_\C \cap F^p E = \bigoplus_{s\ge p} B^{s,m-s} \text{ and }
\pi(F^p E) = \bigoplus_{s\ge p} A^{s,n-s}.
$$
\end{enumerate}

When $A_\Z$ is torsion free, the extension (\ref{ext}) determines an
element $\psi$ of the complex torus $J(\Hom(A,B))$. This is
done as follows: by the property of $\pi$, there is a section
$s_F:A_\C \to E_\C$ that preservers the Hodge filtration; since $A_\Z$ is
torsion free, there is an integral section $s_\Z : A_\Z \to E_\Z$ of $\pi$.
The coset $\psi$ of $s_F - s_\Z$ in $J(\Hom(A,B))$ is independent of the
choices $s_F$ and $s_\Z$.

\subsection{The Theorem of Carlson-Clemens-Morgan}
This is the first example in which periods of (non-abelian) homotopy groups
were related to algebraic cycles.

Here $X$ is a simply connected projective manifold. By
Corollaries~\ref{ex_pi3} and~\ref{mhs_exact}, the sequence
\begin{equation}
\label{pi3_ext}
0 \to H^3(X;\Z)/(\text{torsion}) \to \Hom(\pi_3(X),\Z) \to K \to 0
\end{equation}
is an extension of $\Z$-mixed Hodge structures,\footnote{The integral
statement is proved in \cite{ccm} --- however, $H^3(X;\Z)$ is implicitly
assumed to be torsion free.} where where $K$ is the
kernel of the cup product
$$
S^2 H^2(X;\Z) \to H^4(X;\Z).
$$

Denote the class of a divisor $D$ in the Neron-Severi group
$$
NS(X) := \{\text{group of divisors in $X$}\}/(\text{homological equivalence})
$$
of $X$ by $[D]$. If the codimension 2 cycle
$$
Z:= \sum_{j,k} n_{jk}\, D_j\cap D_k
$$
is homologically trivial, where the $n_{jk}$ are integers and the $D_j$
divisors, then
$$
\widehat{Z} := \sum_{j,k} n_{jk}\, [D_j][D_k] \in S^2 H^2(X;\Z)
$$
is an integral Hodge class of type $(2,2)$ in $K$. Pulling back the extension
(\ref{pi3_ext}) along the mapping $\Z(-2) \to K$ that takes 1 to $\widehat{Z}$,
we obtain an extension
$$
0 \to H^3(X;\Z(2)) \to E_Z \to \Z \to  0
$$
of mixed Hodge structures. This determines a point
$$
\phi_Z \in J(H^3(X;\Z(2))) = J^2(X).
$$

On the other hand, the homologically trivial cycle $Z$ determines a point
$$
\Phi(Z) \in J^2(X).
$$

\begin{theorem}[Carlson-Clemens-Morgan]
The points $\phi_Z$ and $\Phi_Z$ of $J^2(X)$ are equal.
\end{theorem}

\subsection{The Harris-Pulte Theorem}
\label{harris-pulte}
Pulte \cite{pulte} reworked Harris' work on harmonic volume using the Hodge
theory of $\delbar$ and the language of mixed Hodge theory.

Suppose that $C$ is a compact Riemann surface and that $x\in X$.
Corollaries~\ref{ex_pi1} and \ref{mhs_exact} imply that the sequence
$$
0 \to H^1(C) \to H^0(I\Ch_2^\dot(P_{x,x}C)) \to K \to 0
$$
is exact in the category of $\Z$-mixed Hodge structures, where $K$ is the
kernel of the cup product $H^1(C)\otimes H^1(C) \to H^2(C)$. It therefore
determines an element $m_x$ of
$$
J(\Hom(K,H^1(C))).
$$
An element of $\Hom(K,H^1(C))$ can be computed using the recipe in the
previous paragraph. For example, if 
$$
u := \sum_{j,k} a_{jk}\, [w_j]\otimes[\wbar_k] \in K
$$
where each $w_j$ is holomorphic, then, by Proposition~\ref{strictness}, there
is $\xi \in F^1E^1(C)$ such that
$$
d\xi + \sum_{j,k} a_{jk}\, w_j\wedge\wbar_k = 0.
$$
Thus
$$
\int \big(\sum_{j,k} a_{jk}\, w_j\wbar_k + \xi\big)
\in F^1H^0(I\Ch_2^0(P_{x,x}C)).
$$
($s_F$ can be chosen so that this is $s_F(u)$.) The value of the extension
class $\psi$ on $u$ is represented by the homomorphism $H_1(C) \to \C$ obtained
by evaluating this integral on loops based at $x$ representing a basis of
$H_1(C;\Z)$. (Full details can be found in \cite{hain:bowdoin}.)
These integrals are examples of the $\delbar$ analogues of those considered
by Harris.

On the other hand, one has the algebraic 1-cycles
$$
C_x := \big\{[z] - [x] : z \in C\big\} \text{ and }
C_x^- := \big\{[x] - [z] : z \in C\big\}
$$
in the jacobian $\Jac C$ of $C$. These share the same homology class, so the
algebraic cycle
$$
Z_x := C_x - C_x^-
$$
is homologically trivial and determines a point
$$
\nu_x \in J_1(\Jac C) = J(\Lambda^3 H_1(C)(-1)).
$$
The linear mapping $\Lambda^3 H_1(C) \to K^\ast \otimes H_1(C)$
defined by
$$
a\wedge b \wedge c \mapsto \big\{u \mapsto \int_{a\times b}u\big\}\otimes c
+ \big\{u \mapsto \int_{b\times c}u\big\}\otimes a
+ \big\{u \mapsto \int_{c\times a}u\big\}\otimes b
$$
is an injective morphism of Hodge structures, and induces an injection
$$
A : J_1(\Jac C) \hookrightarrow J(\Hom(K,H^1(C))).
$$

\begin{theorem}[Harris-Pulte \cite{harris,pulte}]
With notation as above, $\nu_x = 2A(m_x)$.
\end{theorem}

\begin{remark}
If $C$ is hyperelliptic and $x$ and $y$ are two distinct Weierstrass points,
the mixed Hodge structure on $J(C-\{y\},x)/J^3$ is of order 2. In this case
Colombo \cite{colombo} constructs an extension of $\Z$ by the primitive part
$PH_2(\Jac C;\Z)$ of $H_2(\Jac C)$ from the MHS on $J(C-\{y\},x)/J^4$ and
shows that it is the class of the Collino cycle \cite{collino}, an element
of the Bloch higher Chow group $CH^g(\Jac C,1)$. This example shows that the
MHS on $\pi_1(C-\{y\},x)$ of a hyperelliptic curve contains information about
the extensions associated to elements of higher $K$-groups, ($K_1$ in this
case), not just $K_0$.
\end{remark}

\section{Green's Observation and Conjecture}

Mark Green (unpublished) has given an interpretation of the
Carlson-Clemens-Morgan~Theorem. He also suggested a general picture relating
the Hodge theory of homotopy groups to intersections of cycles. In this
section, we briefly describe Green's ideas, then state and sketch a proof
of a modified version.

\subsection{Green's interpretation}
If one wants to understand the product
$$
CH^a(X) \otimes CH^b(X) \to CH^{a+b}(X)
$$
the first thing one may look at is:
$$
CH^a(X) \otimes CH^b(X) \to \G H^{2a+2b}(X;\Z(a+b))
$$
After this, one may consider:
\begin{multline}
\label{crossover}
\ker\big\{CH^a(X) \otimes CH^b(X) \to \G H^{2a+2b}(X;\Z(a+b)) \big\} \cr
\to J^{a+b}(X) = \Ext^1_\H(\Z,H^{2a+2b-1}(X;\Z(a+b))).
\end{multline}

What Green observed is that when $X$ is a simply connected projective manifold
and $a=b=1$, the result of Carlson-Clemens-Morgan implies this mapping is
determined by the class
$$
\e(X) \in \Ext^1_\H(K,H^3(X;\Z(2)))
$$
of the extension
$$
0 \to H^3(X,\Z(2)) \to \Hom(\pi_3(X),\Z(2)) \to K \to 0,
$$
where $K$ is the kernel of the cup product $S^2H^2(X,\Z(1)) \to H^4(X,\Z(2))$.
This works as follows: since the diagram
$$
\begin{CD}
CH^1(X)\otimes CH^1(X) @>>> H^{2}(X;\Z(1))\otimes H^{2}(X;\Z(1)) \cr
@VVV	@VVV \cr
CH^{2}(X) @>>> H^{4}(X;\Z(2))
\end{CD}
$$
commutes, there is a natural mapping
$$
\ker\big\{CH^1(X) \otimes CH^1(X) \to \G H^{4}(X,\Z(2)) \big\}
\to \G K(2).
$$
The result of Carlson-Clemens-Morgan implies that cupping this homomorphism
with $e(X)$ gives the mapping (\ref{crossover}).

He went on to conjecture that all  the ``crossover mappings'' (\ref{crossover})
--- more generally, all crossover mappings associated to the standard
conjectured filtration of the Chow groups of $X$ ---
are similarly described by cupping with extensions one obtains from the mixed
Hodge structure on homotopy groups of $X$. In his thesis \cite{archava},
Archava proves that a conjecture of Green and Griffiths implies the analogue
of Green's conjecture in the case where the category of mixed Hodge structures
is replaced by the category of arithmetic Hodge structures of Green and
Griffiths \cite{green-griffiths}.

\subsection{Iterated integrals and crossover mappings} This section
proposes a generalization of the theorem of Carlson-Clemens-Morgan to
cycles of all codimensions and also to algebraic manifolds which may be
neither compact nor simply connected.

Suppose that $X$ is a complex algebraic manifold. By Corollary~\ref{mhs_exact},
the sequence
\begin{multline}
\label{ext_mhs}
0 \to QH^{2d-1}(X) \to H^{2d-2}(I\Ch_2^\dot(P_{x,x}X)) \to
[H^{>0}(X)^{\otimes 2}]^{2d}\cr
\stackrel{\text{cup}}{\longrightarrow} H^{2d}(X) \to QH^{2d}(X) \to 0.
\end{multline}
is exact in the category of $\Z$-mixed Hodge structures. Denote by
$H^\ev(X;\Z)$ the sum of the even integral cohomology groups of $X$ of
positive degree. Let
$$
K^\ev = \ker\{H^\ev(X;\Z)^{\otimes 2} \to H^\ev(X;\Z)\}.
$$
This underlies a graded $\Z$-Hodge structure. 
We can pull the extension (\ref{ext_mhs}) back along $K^\ev \to K$ to obtain
a new extension
\begin{equation}
\label{pruned}
0 \to QH^{2d-1}(X;\Z) \to E \to K^\ev \to 0
\end{equation}
which underlies an extension of MHS, which can be seen to be independent of
the choice of the basepoint $x$. There is a natural mapping
$$
\ker\big\{\sum_{\substack{a+b=d \cr a,b>0}}
CH^a(X) \otimes CH^b(X) \to \G H^{2d}(X,\Z(d)) \big\}
\to \G K^{2d}(d).
$$
This, the quotient mapping $H^\dot(X) \to QH^\dot(X)$, and the extension
(\ref{pruned}) determine a homomorphism
\begin{multline*}
\Phi : \ker\big\{\sum_{\substack{a+b=d \cr a,b>0}}
CH^a(X) \otimes CH^b(X) \to H^{2d}(X,\Z(d))\big\} \cr
\to \Ext_\H^1(\Z,QH^{2d-1}(X;\Z(d))).
\end{multline*}

The following, if proven, will generalizes the theorem of Carlson, Clemens
and Morgan.

\begin{conjecture}
If $X$ is a quasi-projective complex algebraic manifold,
the mapping $\Phi$ equals the composition of the crossover mapping
(\ref{crossover}) with the quotient mapping $J^d(X) \to J(QH^{2d+1}(X)(d))$.
\end{conjecture}

\begin{proof}[Heuristic Argument]
By resolution of singularities, we may suppose that the
quasi-projective algebraic manifold $X$ is of the form $\Xbar -D$, where
$\Xbar$ is a complex projective manifold and $D$ is a normal crossings
divisor. Suppose that $Z_1,\dots, Z_m$ are proper algebraic subvarieties of
$X$ of positive codimensions $c_1,\dots,c_m$, respectively. By the moving
lemma, we may move them within their rational equivalence classes
so that they all meet properly. Suppose that the $n_{jk}$ are integers
and that the cycle
$$
W = \sum_{j,k} n_{jk} Z_j \cdot Z_k
$$
is homologically trivial in $X$ of pure codimension $d$.

The basic idea of the argument is easy. The extension class associated to
$W$ is the difference $s_F(\W)-s_\Z(\W)$ mod $F^d$ of Hodge and
integral lifts of the class
$$
\W := \sum_{j,k} n_{jk} [Z_j] \otimes [Z_k] \in K^{2d}
\text{ to }
W_{2d} H^{2d-2}(I\Ch_2^\dot(P_{x,x}X)).
$$
Suppose that $w_j \in E^{c_j,c_j}(\Xbar)$ is a smooth form representing the
Poincar\'e dual of the closure of $Z_j$ in $\Xbar$. Since $W$ is homologically
trivial, there is a form $\xi \in F^d W_1 E^{2d-1}(\Xbar\log D)$ satisfying
$d\xi = \sum n_{jk}\, w_j \wedge w_k$ (cf.\ \cite[I.3.2.8]{hain:dht}.%
\footnote{The $\le$ there should be an equals.}) It follows that
$$
\sum_{j,k} n_{jk}\int w_j w_k + \int \xi
\in F^d W_{2d} H^{2d-2}(I\Ch_2^\dot(P_{x,x}X)),
$$
which we take as the Hodge lift $s_F(\W)$ of $\W$.

In the integral version, we shall use King's theory of logarithmic currents
\cite{king,king:log}.
We would like to take the integral lift of $\W$ to be
\begin{equation}
\label{expr}
s_\Z(\W) := \int \sum_{j,k} n_{jk}\, \d_j \d_k - \d_\G \in
W_{2d} H^{2d-2}(I\Ch_2^\dot(P_{x,x}X))_\Z,
\end{equation}
where $\G$ is a chain of codimension $2d-1$ whose boundary is $W$, and
$\d_j$ is the integration current defined by $Z_j$. To make this argument
precise, one has to show that $s_\Z(\W)$ makes sense. Assume this.

The final task is to compute the extension data. Denote the complex of
currents on $\Xbar$ by $D^\dot(\Xbar)$, and King's complex of
log currents for $(\Xbar,D)$ by $D^\dot(\Xbar\log D)$. These have natural 
Hodge and weight filtrations. There is a log current
$$
\psi_j \in F^{c_j}W_1 D^{2c_j}(\Xbar\log D)
$$
such that $d\psi_j = w_j - \d_j$. Using the formula for the differential,
we have:
\begin{align*}
& \quad \int \big(\w_j w_k - \d_j \d_k \big) \cr
&= -d\int \big(\psi_j \d_k - \d_j \psi_k + \psi_j d\psi_k \big)
- \int
\big(\psi_j \wedge \d_k +  \d_j \wedge \psi_k + \psi_j \wedge d\psi_k \big)\cr
&\equiv 
- \int \big(
\psi_j \wedge \d_k + \d_j \wedge \psi_k + \psi_j \wedge d\psi_k
\big)
\mod \text{ exact forms}
\end{align*}
Combing this with the relations
$$
d\xi = \sum_{j,k} n_{jk}\, w_j \wedge w_k \text{ and }
d\d_\G = - \sum_{j,k} n_{jk}\, \d_j \wedge \d_k
$$
we have, modulo exact forms,
\begin{align*}
s_F(\W) - s_\Z(\W) &\equiv \int \big(\xi + \d_\G \big)
- \sum_{j,k} n_{jk} \int \big(\psi_j\wedge \d_k + \d_j\wedge \psi_k+ \psi_j
\wedge d\psi_k\big) \cr
&\equiv \int_\G \mod F^d + \text{ exact forms,}
\end{align*}
which is the desired result.
\end{proof}

The deficiency in this argument is that the theory of iterated integrals of
currents is not rigorous. To make this argument rigorous, it would be
sufficient to show
that there is a complex of chains whose elements are transverse to
$s_\Z(\W)$, on which $s_\Z(\W)$ takes integral values, and that computes
the integral structure on $H^\dot(I\Ch_2^\dot(P_{x,x}X))$. One possible way to
approach this is to triangulate $\Xbar$ so that $D$, each $Z_j$ and $\G$
are subcomplexes, and then to obtain the cycles that give the integral
structure from some analogue of Adams-Hilton construction \cite{adams-hilton}
associated to the dual cell decomposition. So far, I have not been able to
make this work.

This argument suggests that it is the Hodge theory of iterated integrals
(or more generally, the cosimplicial cobar construction) rather than homotopy
groups which determines periods associated to algebraic cycles, as this
result holds even when the loop space de~Rham theorem and rational homotopy
theory fail. It would be interesting to have an example of an acyclic complex
projective manifold where $\Phi$ is non-trivial to illustrate this point.

This argument also applies in the relative case where the variety $X$ and
the cycles are defined and flat over a smooth base $S$. In this
case, the map $\Phi$ will take values in
$$
\Ext_{\H(S)}^1(\Z_S,R^{2d-1}f_\ast\,\Z_X(d))
$$
where $f:X \to S$ and $\H(S)$ denotes the category of admissible variations
of mixed Hodge structure over $S$. This can be seen using results from
\cite[Part~II]{hain:dht} and \cite{hain-zucker}.
By combining this with the standard technique of spreading a variety defined
over a subfield of $\C$, one should get elements of the Hodge realization of
motivic cohomology as considered in \cite{asakura}, for example.

\section{Beyond Nilpotence}
\label{beyond}

The applicability of Chen's de~Rham theory (equivalently, rational homotopy
theory) is limited by nilpotence. Using ordinary iterated line integrals,
one can only separate those elements of $\pi_1(X,x)$ that can be separated
by homomorphisms from $\pi_1(X,x)$ to a group of unipotent upper triangular
matrices. If the first Betti number $b_1(X)$ of $X$ is zero, all such
homomorphisms are trivial, and iterated line integrals cannot separate any
elements of $\pi_1(X,x)$ from the identity. If $b_1(X) = 1$, then
the image of all such homomorphisms is abelian, and iterated line integrals
can separate only those elements that are distinct in $H_1(X;\R)$. Thus,
in order to apply de~Rham theory to the study of moduli spaces of curves
and mapping class groups ($b_1(X)=0$) or knot groups ($b_1(X)=1$), for 
example, iterated integrals need to be generalized.

Before explaining two ways of doing this we shall restate Chen's de~Rham
theorem for the fundamental group in a form suitable for generalization.

First recall the definition of unipotent (also known as Malcev) completion.
A unipotent group is a Lie group that can be realized as a closed
subgroup of the group of a group of unipotent upper triangular matrices.
(That is, upper triangular matrices with 1's on the diagonal.) Unipotent
groups are necessarily algebraic groups as the exponential map from the
Lie algebra of strictly upper triangular matrices to the group of unipotent
upper triangular matrices is a polynomial bijection.%
\footnote{Here and below, I shall be vague about the field $F$ of definition
of the group. It will always be either $\R$ or $\C$. Also, I will not
distinguish between the algebraic group and its group of $F$-rational points.}

Suppose that $\G$ is a discrete group. A homomorphism $\rho$ from $\G$ to a 
unipotent group $U$ is Zariski dense if there is no proper unipotent subgroup
of $U$ that contains the image of $\rho$. The set of Zariski dense unipotent
representations $\rho : \G \to U_\rho$ forms an inverse system. The {\it
unipotent completion} of $\G$ is the inverse limit of all such representations;
it is a homomorphism from $\G$ into the {\it prounipotent} group
$$
\U(\G) := \limproj \rho U_\rho.
$$
Every homomorphism $\G \to U$ from $\G$ to a unipotent group factors
through the natural homomorpism $\G \to \U(\G)$.
The coordinate ring of $\U(\G)$ is, by definition, the direct limit of the
coordinate rings of the $U_\rho$:
$$
\O(\U(\G)) = \liminj \rho \O(U_\rho).
$$
It is isomorphic to the Hopf algebra of matrix entries $f : \G \to \R$ of
all unipotent representations of $\G$.

The following statement is equivalent to the statement of Chen's de~Rham
theorem for the fundamental group given in Section~\ref{derham}.

\begin{theorem}
If $X$ is a connected manifold, then integration induces a Hopf algebra
isomorphism
$$
\O(\U(\pi_1(X,x))) \cong H^0(\Ch^\dot(P_{x,x}X)).
$$
\end{theorem}

One recovers the unipotent completion of $\pi_1(X,x)$ as
$\mathrm{Spec}H^0(\Ch^\dot(P_{x,x}X))$. The homomorphism
$\pi_1(X,x) \to \U(\pi_1(X,x))$ takes the homotopy class of the loop $\gamma$
to the maximal ideal of iterated integrals that vanish on it.

\subsection{Relative unipotent completion}
\label{rel_comp}
Deligne suggested the following generalization of unipotent completion, which
is itself a generalization of the idea of the algebraic envelope of a 
discrete group defined by Hochschild and Mostow \cite[\S4]{hoch-mostow}.

Suppose that $S$ is a reductive algebraic group. (That is, an affine
algebraic group, all of whose finite dimensional representations are
completely reducible, such as $SL_n$, $GL_n$, $O(n)$, $\Gm$, \dots .)
Suppose that $\G$ is a discrete group as above and that $\rho : \G \to S$
is a Zariski dense homomorphism.

Similar to the construction of the unipotent completion of $\G$, one
can construct a proalgebraic group $\cG(\G,\rho)$, which is an extension
$$
1 \to \U(\G,\rho) \to \cG(\G,\rho) \stackrel{p}{\to} S \to 1
$$
of $S$ by a prounipotent group, and a homomorphism $\G \to \cG(\G,\rho)$
whose composition with $p$ is $\rho$. Every homomorphism from $\G$ into an
algebraic group $G$ that is an extension of $S$ by a unipotent group, and
for which the composite $\G \to G \to S$ is $\rho$, factors through
the natural homomorphism $\G \to \cG(\G,\rho)$.

The homomorphism $\G \to \cG(\G,\rho)$ is called the {\it completion of $\G$
relative to $\rho$}.  When $S$ is trivial, the relative completion reduces
to classical unipotent completion described above.

The definition of iterated integrals can be generalized to more general
forms to compute the coordinate rings of relative completions of fundamental
groups.
Suppose now that $\G = \pi_1(X,x)$, where $X$ is a connected manifold.
The representation $\rho$ determines a flat principal $S$-bundle, $P \to X$,
together with an identification of the fiber over $x$ with $S$. One can then
consider the corresponding (infinite dimensional) bundle $\O(P) \to X$
whose fiber over $y \in X$ is the coordinate ring of the fiber of $P$
over $y$. This is a flat bundle of $\R$-algebras. One can, consider the
dga $E^\dot(X,\O(P))$ of $S$-finite differential
forms on $X$ with coefficients in $\O(P)$. In \cite{hain:malcev}, Chen's
definition of iterated integrals is extended to such forms. The iterated
integrals of degree 0 are, as before, functions $P_{x,x}X \to \R$.

Two augmentations
$$
\delta : E^\dot(X,\O(P)) \to \O(S) \text{ and }
\epsilon : E^\dot(X,\O(P)) \to \R
$$
are obtained by restricting forms to the fiber $S$ over $x$ and to the
identity $1\in S$ in this fiber. These, give $\O(S)$ and $\R$ structures
of modules over $E^\dot(X,\O(P))$. One can then form the bar construction
$$
B(\R,E^\dot(X,\O(P)),\O(S)).
$$
This maps to the complex of iterated integrals of elements
of $E^\dot(X,\O(P))$.

\begin{theorem}
Integration of iterated integrals induces a natural isomorphism
$$
H^0(B(\R,E^\dot(X,\O(P)),\O(S))) \cong \O(\cG(\pi_1(X,x),\rho)).
$$
\end{theorem}

The corresponding Hodge theory is developed in \cite{hain:malcev}. It
is used in \cite{hain:torelli} to give an explicit presentation of the
completion of mapping class groups $\G_g$ with respect to the standard
homomorphism $\G_g \to Sp_g$ to the symplectic group given by the action
of $\G_g$ on the first homology of a genus $g$ surface when $g\ge 6$.

One disadvantage of the generalization sketched above is that these
generalized iterated integrals, being constructed from differential forms
with values in a flat vector bundle, are not so easy to work with. A more
direct and concrete approach is possible in the solvable case.

\subsection{Solvable iterated integrals}
In his senior thesis, Carl Miller \cite{miller} considers the solvable case.
Here it is best to take the ground field to be $\C$. The reductive group
is a diagonalizable algebraic group:
$$
S=(\C^\ast)^k \times \mu_{d_1}\times \cdots \times \mu_{d_m}.
$$
He defines {\it exponential iterated line integrals}, which are certain
convergent infinite sums of standard iterated line integrals of the type
that occur as matrix entries of solvable representations of fundamental
groups. Exponential iterated line integrals are linear combinations of
iterated line integrals of the form
\begin{multline*}
\int e^{\d_0} w_1 e^{\d_1} w_2 e^{\d_3} \dots e^{\d_{n-1}} w_n e^{\d_n} \cr
:= \sum_{k_j \ge 0} \int
\overbrace{\d_0 \dots \d_0}^{k_0}
w_1
\overbrace{d_1 \dots d_1}^{k_1}
w_2
\overbrace{\d_2 \dots \d_2}^{k_2}
\dots
\overbrace{\d_{n-1} \dots \d_{n-1}}^{k_{n-1}}
w_n
\overbrace{\d_n \dots \d_n}^{k_n}
\end{multline*}
where $\d_0,\dots, \d_n, w_1, \dots, w_n$ are all 1-forms. This sum
converges absolutely when evaluated on any path. The terminology and notation
derive from the easily verified fact that
$$
\exp \int_\gamma w = \sum_{k\ge 0}\, \int_\gamma \overbrace{w\dots w}^k.
$$

\begin{theorem}[Miller]
Suppose that $X$ is a connected manifold and
$$
\rho : \pi_1(X,x) \to S
$$
is a Zariski dense representation to a diagonalizable $\C$-algebraic group.
If $\rho$ factors through $H_1(X)/\text{torsion}$, then
the Hopf algebra of locally constant exponential iterated integrals
associated to $\rho$ is isomorphic to the coordinate ring
$\O(\cG(\pi_1(X,x),\rho))$ of the completion of $\pi_1(X,x)$ relative to
$\rho$.
\end{theorem}

He also shows that for a large class of knots $K$ (which includes all fibered
knots), there is a representation $\rho:\pi_1(S^3-K,x) \to S$ into a
diagonalizable algebraic group such that $\pi_1(S^3-K,x)$ injects into the
corresponding relative completion. In particular, there are enough exponential
iterated line integrals to separate elements of $\pi_1(S^3-K,x)$. The
representation $\rho$ can be computed from the eigenvalues of the Alexander
polynomial of $K$. The representation $\pi_1(S^3-K,x) \to S$ is the Zariski
closure of the ``semi-simplification'' of the Alexander module of $K$.

\section{Algebraic Iterated Integrals}
\label{alg_dr}

A standard tool in the study of algebraic varieties over any field is
algebraic de~Rham theory, which originates in the theory of Riemann
surfaces and was generalized by Grothendieck \cite{grothendieck} among 
others. This algebraic de~Rham theory extends to iterated integrals and
several approaches will be presented in this section. I will begin with the
most elementary and progress to the abstract, but powerful, approach of
Wojtkowiak \cite{wojtkowiak}.

\subsection{Iterated integrals of the second kind}

The historical roots of algebraic de~Rham cohomology lie in the classical
result regarding differentials of the second kind on a compact Riemann
surface. Recall that a meromorphic 1-form $w$ on a compact Riemann surface
$X$ is of the {\it second kind} if it has zero residue at each point.
Alternatively, $w$ is of the second kind if the value of $\int w$ on each
loop in $X-\{\text{singularities of } w\}$ depends only on the class of
the loop in $H_1(X)$. A classical result asserts that there is a natural
isomorphism
$$
H^1(X;\C) \cong\,
\frac{\{\text{meromorphic differentials of the second kind on $X$}\}}
{\{\text{differentials of meromorphic functions}\}}
$$

This can be generalized to iterated integrals. Suppose that $X$ is a compact
Riemann and that $S$ is a finite subset. An {\it iterated line integral
of the second kind} on $X-S$ is an iterated integral
$$
\sum_{r \le s} \sum_{|I| = r} a_I \int w_{i_i} \dots w_{i_r},
$$
where $a_I\in \C$ and each $w_j$ is a meromorphic differential on $X$, with
the property that its value on each path in $X$ that avoids the singularities
of all $w_j$ depends only on its homotopy class (relative to its endpoints)
in $X-S$.

\begin{example}[{cf.\ \cite[p.~260]{hain:bowdoin}}]
We will assume that $S$ is empty (the case where $S$ is non empty is simpler).
Suppose that $w_1,\dots, w_n$ are differentials of the second kind on $X$ and
that $a_{jk} \in \C$. Since differentials of the second kind are locally (in
the complex topology) the exterior derivative of a meromorphic function, for
each point $x \in U$ we can find a function $f_j$, meromorphic at $x$, such
that $d f_j = w_j$ about $x$. Define
$$
r_{jk}(x) = \Res_{z=x} \big[f_j(z)w_k(z)\big].
$$
Since $w_k$ is of the second kind, changing $f_j$ by a constant will not
change $r_{jk}(x)$. If
$$
\sum_{x \in X} \sum_{j,k} a_{jk}r_{jk}(x) = 0
$$
there is a meromorphic differential $u$ on $X$ (which can be taken to be
of the third kind) such that
$$
\Res_{z=x} u(z) = - \sum_{j,k} a_{jk}r_{jk}(x).
$$
The iterated integral
$$
\sum_{j,k} a_{jk} \int w_j w_k + \int u
$$
is of the second kind. This can be seen by noting that the integrand
of this integral near $x$ is
$$
\sum_{j,k} a_{jk} f_j(z)w_k(z) + u(z),
$$
which has zero residue at $x$. Equivalently, the pullback of the integrand
of this iterated integral to the universal covering of $X-S$ is of the
second kind.
\end{example}

\begin{theorem}
\label{sec_kind}
If $X$ is a compact Riemann surface, $S$ a finite subset of $X$, and
$1\le s \le \infty$, then, for all $x,y \in X-S$, integration induces a
natural isomorphism
$$
H^0(\Ch_s^\dot(P_{x,y}(X-S))_\C) \cong
\left\{
\parbox{2.1in}{The set of iterated integrals of the second kind
of length $\le s$ on $X-S$}
\right\}
$$
\end{theorem}

\begin{proof}
This is just an algebraic version of the proof of Chen's $\pi_1$ de~Rham
theorem given in \cite[\S4]{hain:bowdoin}. Familiarity with that proof will
be assumed. I will just make those additional points necessary to prove this
variant.

Set $U=X-S$. Suppose that $s < \infty$. We consider the truncated group ring
$\C\pi_1(U,x)/J^{s+1}$ to be a $\pi_1(U,x)$-module via right multiplication.
Let $E_s \to U$ be the corresponding flat bundle.
This is a holomorphic vector bundle with a flat holomorphic connection.
It is filtered by the flat subbundles corresponding to filtration
$$
\C\pi_1(U,x)/J^{s+1} \supseteq J/J^{s+1}
\supseteq \dots \supseteq J^s/J^{s+1} \supseteq 0
$$
of $\C\pi_1(U,x)/J^{s+1}$ by right $\pi_1(U,x)$-submodules. Denote the
corresponding filtration of $\E$ by
$$
\E_s = \E_s^0 \supseteq \E_s^1 \supseteq \cdots \supseteq \E_s^s \supseteq 0.
$$
By the calculation in \cite[Prop.~4.2]{hain:bowdoin}, each of the bundles
$\E_s^t/\E_s^{t+1}$ has trivial monodromy, so that each $\E_s^t$ has 
unipotent monodromy.

By the results of \cite{deligne:ode}, each of the bundles $\E_s^t$ has a
canonical extension $\overline\E_s^t$ to $X$. These satisfy:
\begin{enumerate}
\item each $\Ebar_s^t$ is a subbundle of $\Ebar_s:=\Ebar_s^0$;
\item the connection on $\E_s$ extends to a meromorphic connection on $\Ebar_s$
which restricts to a meromorphic connection on each of the $\Ebar_s^t$;
\item the connection on each of the bundles $\Ebar_s^t/\Ebar_s^{t+1}$ is trivial
over $X$.
\end{enumerate}
(Take $\Ebar_s^t = \E_s^t$ when $S$ is empty.)

The following lemma implies that there are meromorphic
trivializations of each $\Ebar_s$ compatible with all of the projections
$$
\cdots \to \Ebar_s \to \Ebar_{s-1} \to \cdots \to \Ebar_0 = \O_X
$$
and where the induced trivializations of each graded quotient of each
$\Ebar_s$ is flat. Moreover, we can arrange for all of the singularities
of the trivialization\footnote{A meromorphic trivialization
$\phi : E \to \O_X^N$ is singular at $x$ if either $\phi$ has a pole at $x$
or if the determinant of $\phi$ vanishes at $x$.} to lie in any
prescribed non-empty finite subset $T$ of $X$.

The connection form $\w_s$ of $\Ebar_s$ with respect to this trivialization
thus satisfies
$$
\w_s \in \{\text{meromorphic 1-forms on $X$}\}
\otimes J^{-1}\End(\C\pi_1(U,x)/J^{s+1})
$$
with values in the linear endomorphisms of $\C\pi_1(U,x)/J^{s+1}$ that
preserve the filtration
$$
\C\pi_1(U,x)/J^{s+1} \supseteq J/J^{s+1}
\supseteq \dots \supseteq J^s/J^{s+1} \supseteq 0
$$
and act trivially on its graded quotients. This connection is thus nilpotent.
Note that, even though $\w_s$ may have  poles in $U$, the connection given
by $\w_s$ has trivial monodromy about each point of $U$. This is the key point
in the proof; it implies that the transport \cite[\S2]{hain:bowdoin}
$$
T =
1 + \int \w_s + \int \w_s \w_s + \cdots + \int \overbrace{\w_s \dots \w_s}^s
$$
is an $\End \C\pi_1(U,x)/J^{s+1}$-valued iterated integral of the second kind
on $U$. Its matrix entries are iterated integrals of the second kind.

The result when $x=y$ now follows as in the proof of \cite[\S4]{hain:bowdoin}.
The case when $x \neq y$ is easily deduced from the case $x=y$.
The result for $s=\infty$ is obtained by taking direct limits over $s$ using
the fact that $\w_s$ is the image of $\w_{s+1}$ under the projection
\begin{multline*}
\{\text{meromorphic 1-forms on $X$}\}
\otimes J^{-1}\End(\C\pi_1(U,x)/J^{s+2}) \cr
\to
\{\text{meromorphic 1-forms on $X$}\}
\otimes J^{-1}\End(\C\pi_1(U,x)/J^{s+1})
\end{multline*}
\end{proof}

\begin{lemma}
Suppose that
$$
0 \to \O_X^N \to \E \stackrel{p}{\to} \F \to 0
$$
is an extension of holomorphic vector bundles over a compact Riemann surface
$X$. If $T$ is a non-empty subset of $X$, there is a meromorphic splitting
of $p$ which is holomorphic outside $T$.
\end{lemma}

\begin{proof}
Set $\Fdual = \Hom(\F,\O_X)$.
Riemann-Roch implies that $H^1(X,\Fdual(\ast T)) = 0$, where $\Fdual(\ast T)$
is defined to be the sheaf of meromorphic sections of $\Fdual$ that are
holomorphic outside $T$. It follows from obstruction theory for extensions of
vector bundles that the sequence has a meromorphic splitting that is holomorphic
on $X-T$.
\end{proof}

\begin{remark}
Note that if $S$ is non-empty, the proof shows that the algebraic iterated
line integrals built out of meromorphic forms that are holomorphic on 
$X-S$ equals $H^0(\Ch^\dot(P_{x,y}(X-S))_\C)$. Since $X-S$ is affine, this
is a very special case of Theorem~\ref{groth_dr} in the next paragraph,
a consequence of Grothendieck's algebraic de~Rham Theorem. The result
above can also be used to show that if $X$ is a smooth curve defined over
a subfield $F$ of $\C$, then $H^0(\Ch^\dot(P_{x,y}(X-S))_\C)$ has a 
canonical $F$-form --- namely that consisting of those meromorphic
differentials of the second kind on $X-S$ that are defined over $F$.
\end{remark}

It would be interesting and useful to have a description of the Hodge and
weight filtrations on $H^0(\Ch^\dot(P_{x,y}(X-S))_\C)$, possibly in terms
of some kind of pole filtration, as one has for cohomology.

\subsection{Grothendieck's theorem and its analogues for iterated integrals}

Suppose that $X$ is a variety over a field $F$ of characteristic zero.
Denote the sheaf of K\"ahler differentials of $X$ over $F$ by
$\Omega_{X/F}^\dot$. Denote its global sections over $X$ by
$H^0(\Omega^\dot_{X/F})$. It is a commutative differential graded algebra
over $F$. When $F=\C$ and $X$ is smooth, every algebraic differential
$w \in H^0(\Omega^\dot_{X/\C})$ is a holomorphic differential on $X$. The
corresponding mapping $H^0(\Omega^\dot_{X/\C}) \to E^\dot(X)_\C$ is a dga
homomorphism.

\begin{theorem}
\label{groth_dr}
If $X$ is a complex affine manifold, then natural homomorphism
$$
H^\dot(H^0(\Omega^\dot_{X/\C})) \to H^\dot(X;\C)
$$
is a ring isomorphism.
\end{theorem}

Note that Theorem~\ref{sec_kind} is a consequence of this when $S$ is
non-empty and $S=T$.

If $F \subset\C$ and $X$ is defined over $F$, then
$H^0(\Omega^\dot_{X/F})\otimes_F \C \cong H^0(\Omega^\dot_{X/\C})$.
One important consequence of Grothendieck's theorem is that if $F$ is a 
subfield of $\C$, then
$$
H^\dot(X(\C);\C) \cong H^\dot(H^0(\Omega^\dot_{X/F}))\otimes_F \C.
$$
That is, the de~Rham cohomology of the complex manifold $X(\C)$ has a
natural $F$ structure which is functorial with respect to morphisms of
affine manifolds over $F$.

This can be generalized to arbitrary smooth varieties over $F$ by taking
hypercohomology. Define the {\it algebraic de~Rham cohomology} of $X$ by
$$
\HDR^\dot(X) = \bH^\dot(X,\Omega_{X/F}^\dot).
$$
As above, if $F$ is a subfield of $\C$, then the ordinary de~Rham cohomology
of $X(\C)$ has a natural $F$-structure:
$$
H^\dot(X(\C);\C) \cong \HDR^\dot(X)\otimes_F \C.
$$
Using the classical Hodge theorem, one can show that if $X$ is also projective,
the Hodge filtration
$$
F^p H^m(X(\C)) := \bigoplus_{s \ge p} H^{s,m-s}(X(\C))
$$
is obtained from a natural Hodge filtration
$$
\HDR^\dot(X) = F^0\HDR^\dot(X) \supseteq F^1\HDR^\dot(X)
\supseteq F^2 \HDR^\dot(X) \supseteq \cdots
$$
of the algebraic de~Rham cohomology by tensoring by $\C$.

This can be extended to iterated integrals on affine manifolds in the
obvious way. For an affine manifold $X$ over $F$ and $F$-rational points
$x,y \in X(F)$, define the algebraic iterated integrals on $P_{x,y}X$ by
$$
\HDR^\dot(P_{x,y}X) = H^\dot(B(F,H^0(\Omega^\dot_{X/F}),F))
$$
where $F$ is viewed as a module over $H^0(\Omega^\dot_{X/F})$ via the
two augmentations induced by $x$ and $y$.\footnote{This is not an
unreasonable definition, but one should recall that that when $X$ is not
simply connected and $F=\C$, the de~Rham theorem may not hold as we have
seen in Example~\ref{limits}.} It follows from Corollary~\ref{A_ints} and
Grothendieck's theorem above that if $F$ is a subfield of $\C$, there is a
canonical isomorphism
\begin{equation}
\label{defn}
\HDR^\dot(P_{x,y}X)\otimes_F \C \cong H^\dot(\Ch^\dot(P_{x,y}X(\C))).
\end{equation}

When $X$ is not affine, one can replace $X$ by a smooth affine hypercovering
$U_\dot \to X$ and apply the methods of \cite[\S5]{hain:dht} or Navarro
\cite{navarro} (see below) to construct a commutative dga $A^\dot(U_\dot)$
over $F$ with the property that when tensored with $\C$ over $F$, it is
naturally quasi-isomorphic to $E^\dot(X)_\C$. One can then define
$\HDR^\dot(P_{x,y}X)$ to be the cohomology of the corresponding
bar construction as above. It will give a natural $F$-form of
$H^\dot(\Ch^\dot(P_{x,y}X(\C)))$. However, in this general case, it
is better to use Wojtkowiak's approach, which is explained in the next
paragraph.

\subsection{Wojtkowiak's approach}

The most functorial way to approach algebraic de~Rham theory of iterated
integrals on varieties
is via the works of Navarro \cite{navarro} and Wojtkowiak \cite{wojtkowiak}.
This approach has been used in the works of Shiho \cite{shiho} and Kim-Hain
\cite{kim-hain} on the crystalline version of unipotent completion.

Suppose that $D$ is a normal crossing divisor in a smooth complete variety
$\Xbar$, both defined over a field $F$ of characteristic zero. Set $X=\Xbar-D$
and denote the inclusion $X\hookrightarrow \Xbar$ by $j$. One then
has the sheaf of logarithmic differentials $\Omega_\Xbar^\dot(\log D)$ on $\Xbar$,
which is quasi-isomorphic to $j_\ast F$.

For a continuous map $f:U\to V$ between topological spaces,
Navarro \cite{navarro} has constructed a
functor $\RTW f_\ast$ from the category of complexes of sheaves on $U$ to the
category of complexes of sheaves on $V$ with many wonderful properties.
Among them:
\begin{enumerate}
\item $\RTW f_\ast$ takes sheaves of commutative dgas on $U$ to sheaves of
commutative dgas on $V$;
\item if $V$ is a point, then the global sections $\G \RTW f_\ast \Q_U$ of
$\RTW f_\ast \Q_U$ is Sullivan's rational de~Rham complex of $U$;
\item $\RTW f_\ast$ takes quasi-isomorphisms to quasi-isomorphisms;
\item it induces the usual $Rf_\ast$ from the derived category of
bounded complex of sheaves on $U$ to the bounded derived category of sheaves
on $V$.
\end{enumerate}
For convenience, we denote the global sections $\G \RTW$ of $\RTW$ by $\rtw$.
For an arbitrary topological space $Z$, define
$$
A^\dot(Z) = \rtw \Q_Z.
$$
This is the Thom-Whitney-Sullivan de~Rham complex of $Z$. Its cohomology is
naturally isomorphic to $H^\dot(Z;\Q)$.

In the present situation, we can assign the commutative differential graded
algebra
$$
L^\dot(\Xbar,D) := \rtw \Omega^\dot_\Xbar(\log D)
$$
to $(\Xbar,D)$, where we are viewing $\Xbar$ as a topological space in the
Zariski topology. This dga is natural in the pair $(\Xbar,D)$.

If $x,y$ are $F$-rational points of $X$, there are natural augmentations
$L^\dot(\Xbar,D) \to F$. We can therefore use them to form the bar construction
$B(F,L^\dot(\Xbar,D),F)$. Following Wojtkowiak\footnote{Actually, he does not
use logarithmic forms, just algebraic forms on $X$. However, it is necessary
to use logarithmic forms in order to compute the Hodge and weight filtrations.}
\cite{wojtkowiak}, we define
$$
\HDR^\dot(P_{x,y}X) = H^\dot(B(F,L^\dot(\Xbar,D),F)).
$$
This definition agrees with the ones above.

If $F$ is a subfield of $\C$, then the naturality of Navarro's functor
implies that there is a natural dga quasi-isomorphism
$$
A^\dot(X(\C))\otimes_\Q \C \leftrightarrow L^\dot(\Xbar,D)\otimes_F \C
$$
where we regard $X(\C)$ as a topological space in the complex topology.
This quasi-isomorphism respects the augmentations induced by
$x$ and $y$. Thus we have:

\begin{theorem}[Wojtkowiak]
If $F$ is a subfield of $\C$, there is a natural isomorphism
\begin{equation}
\label{isom}
\HDR^\dot(P_{x,y}X)\otimes_F \C \cong H^\dot(\Ch^\dot(P_{x,y}(X(\C)))).
\end{equation}
\end{theorem}

This result can be extended to the Hodge and weight filtrations.
The Hodge filtration of $L^\dot(\Xbar,D)$ is defined by
$$
F^p L^\dot(\Xbar,D) = \rtw
\big[\Omega_\Xbar^p(\log D) \to \Omega_\Xbar^{p+1}(\log D) \to \cdots \big],
$$
where $\Omega_\Xbar^p(\log D)$ is placed in degree $p$.
This extends to a Hodge filtration on $B(F,L^\dot(\Xbar,D),F)$ as described
in \cite[\S3.2]{hain:dht}. The Hodge filtration of $\HDR^\dot(P_{x,y}X)$ is
defined by
\begin{align*}
F^p \HDR^\dot(P_{x,y}X) &= \im\big\{
H^\dot(F^p B(F,L^\dot(\Xbar,D),F)) \to \HDR^\dot(P_{x,y}X)
\big\} \cr
&\cong H^\dot(F^p B(F,L^\dot(\Xbar,D),F)).
\end{align*}

Similarly, the weight filtration of $L^\dot(\Xbar,D)$ is defined by
$$
W_m L^\dot(\Xbar,D) = \rtw \tau_{\le m} \Omega_\Xbar^\dot(\log D).
$$
Like the Hodge filtration, this extends to a weight filtration of
$B(F,L^\dot(\Xbar,D),F)$ as in \cite[\S3.2]{hain:dht}. The weight filtration
of $\HDR^\dot(P_{x,y}X)$ is defined by
$$
W_m \HDR^n(P_{x,y}X) = \im\big\{
H^n(W_{m-n} B(F,L^\dot(\Xbar,D),F)) \to \HDR^\dot(P_{x,y}X)
\big\}.
$$

\begin{theorem}
Suppose, as above, that $\Xbar$ is a smooth complete variety and $D$ a
normal crossings divisor in $\Xbar$, both defined over $F$. If
$X=\Xbar - D$, then there is a Hodge filtration
$$
\HDR^\dot(P_{x,y}X) = F^0 \HDR^\dot(P_{x,y}X)\supseteq \HDR^\dot(P_{x,y}X)
\supseteq \HDR^\dot(P_{x,y}X) \supseteq \cdots
$$
and a weight filtration
$$
\cdots \subseteq W_m\HDR^\dot(P_{x,y}X) \subseteq W_{m+1}\HDR^\dot(P_{x,y}X)
\subseteq \cdots \subseteq HDR^\dot(P_{x,y}X)
$$
which are functorial with respect to morphisms of smooth $F$-varieties and
are compatible with the product and, when $x=y$, the coproduct and antipode.
These filtrations behave well under extension of scalars; that is, if $K$
is an extension field of $F$, then there are natural isomorphisms
$$
F^p \HDR^\dot(P_{x,y}X\otimes_F K) \cong
\big( F^p \HDR^\dot(P_{x,y}X)\big)\otimes_F K
$$
and
$$
W_m \HDR^\dot(P_{x,y}X\otimes_F K) \cong
\big( W_m\HDR^\dot(P_{x,y}X)\big)\otimes_F K.
$$
When $F=\C$, these filtrations agree with those defined in \cite{hain:dht}.
\end{theorem}

\begin{proof}
The first point is that there is a natural filtered
quasi-isomorphism
$$
(E^\dot(\Xbar \log D),F^\dot) \leftrightarrow (L^\dot(\Xbar,D),F^\dot).
$$
The second is that there are natural quasi-isomorphisms
$$
j_\ast F_X \hookrightarrow \Omega_\Xbar^\dot(\log D) \hookrightarrow
j_\ast \Omega_X^\dot.
$$
\end{proof}

\section{The Cobar Construction}
\label{cobar}

In this section, we review the cobar construction (a cosimplicial models of
loop and path spaces) and explain how iterated integrals are the ``de~Rham
realization'' of it. The applications of iterated integrals in earlier
sections, and their role in the algebraic de~Rham theorems for varieties
over arbitrary fields, suggest that the cosimplicial version of the cobar
construction plays a direct and deep role in the theory of motives and that
the examples presented in this paper are just the Hodge-de~Rham realizations
of such motivic phenomena. Additional evidence for this view comes from the
works of Colombo \cite{colombo}, Cushman \cite{cushman}, Shiho \cite{shiho}
and Terasoma \cite{terasoma}.

The original version of the cobar construction, due to Frank Adams
\cite{adams:pna,adams}, grew out of earlier work \cite{adams-hilton} with
Peter Hilton. Adams' cobar construction can be viewed as a functorial
construction which associates to a certain singular chain complex
$\So_\dot(X,x)$ of a pointed space $(X,x)$, a complex $\Ad(\So_\dot(X,x))$
that maps to the reduced {\em cubical} chains on the loopspace $P_{x,x}X$ and
which is dual, in some sense, to the bar construction on the dual of
$\So_\dot(X,x)$. The map from Adams' cobar construction to the reduced cubical
chains is a quasi-isomorphism when $X$ is simply connected. In the non-simply
connected case, a result of Stallings \cite{stallings} implies that
$H_0(\Ad(\So_\dot(X,x)))$ is naturally isomorphic to
$H_0(P_{x,x}X;\Z) = \Z\pi_1(X,x)$.

We begin with the abstract cobar construction and work back
towards the classical one. The abstract approach appears to originate with
the book of Bousfield and Kan \cite{bousfield-kan}. Much of what we write
here is an elaboration of the first section of Wojtkowiak's paper
\cite{wojtkowiak}. Chen has given a nice exposition of the classical cobar
construction in the appendix of \cite{chen:bams}.

\subsection{Simplicial and cosimplicial objects}
Denote the category of finite ordinals by $\simp$; its objects are the
finite ordinals $[n] := \{0,1,\dots,n\}$
and the morphisms are order preserving functions. Among these, the face maps
$$
	d^j : [n-1] \to [n], \quad 0 \le j \le n
$$
play a special role; $d^j$ is the unique order preserving injection that
omits the value $j$.

A contravariant functor $\simp \to \cC$ is called a {\it simplicial object}
in the category $\cC$.  A {\it cosimplicial object} of $\cC$ is a covariant
functor $\simp \to \cC$.

\begin{example}
Denote the standard $n$-simplex by $\D^n$. We can regard its vertices
as being the ordinal $[n]$. Each order preserving mapping $f:[n] \to [m]$
induces a linear mapping $|f| : \D^n \to \D^m$. These assemble to give
the cosimplicial space $\D^\dot$
$$
\xymatrix{
\D^0 \ar@<0.5ex>[r]^{d^0}\ar@<-0.5ex>[r]_{d^1} & \D^1 \ar@<1ex>[r]^{d^0}
\ar[r]\ar@<-1ex>[r]_{d^2}
& \D^2 \ar@<1.5ex>[r]^{d^0} \ar@<0.5ex>[r] \ar@<-0.5ex>[r] \ar@<-1.5ex>[r]_{d^3}
& \D^3 & \cdots
}
$$
whose value on $[n]$ is $\D^n$.
\end{example}

\begin{example}
Suppose that $K$ is an ordered finite simplicial complex (that is, there is
a total order on the vertices of each simplex). Then one has the simplicial
set $K_\dot$ whose set of $n$-simplices $K_n$ is the set of order preserving
mappings $\phi : [n] \to K$ (not necessarily injective) such that the images
of the $\phi(j)$ span a simplex of $K$. In particular, we have the simplicial
set $\Delta^n_\dot$ whose set of $m$-simplices is the set of all order
preserving mappings from $[m]$ to $[n]$.
\end{example}

If one has a simplicial or cosimplicial abelian group, one obtains
a chain complex simply by defining the differential to be the alternating
sum of the (co)face maps. Likewise, if one has a simplicial or cosimplicial
chain complex, one obtains a double complex.

\subsection{Cosimplicial models of path and loop spaces}
\label{cosimp}
Suppose that $X$ is a topological space. Denote the simplicial model of
the unit interval $\D^1_\dot$ by $\I$. Let
$$
X^{\I} = \Hom(\I, X).
$$
This is a cosimplicial space which models the full path space $PX$. Its
space of $n$-cosimplices is $\Hom(I_n,X)$.
Since there are $n+2$ order preserving mappings $[n] \to \{0,1\}$, this is
just $X^{n+2}$. The $j$th coface mapping $d^j : X^{I_{n-1}} \to X^{I_n}$ is
$$
\overbrace{\id\times \cdots \times \id}^{j} \times (\text{diagonal})
\times \overbrace{\id \times \cdots \times \id}^{n-j} : X^{n+1} \to X^{n+2} 
$$
We shall denote it by $P^\dot X$ and its set of $n$-cosimplices by $P^n X$.

The simplicial set $\partial \I$ is the simplicial set associated to the
discrete set $\{0,1\}$. Since $(\partial I)_n$ consists of just the two constant
maps $[n] \to \{0,1\}$, the cosimplicial space $X^{\partial \I}$ consists
of $X\times X$ in each degree. The mapping $X^\I \to X^{\partial \I}$
corresponds to the projection $PX \to X\times X$ that takes a path $\gamma$
to its endpoints $(\gamma(0),\gamma(1))$.

One obtains a cosimplicial model $P^\dot_{x,y}X$ for $P_{x,y}X$ by taking
the fiber of $X^\I \to X^{\partial \I}$. Specifically, $P^n_{x,y} X = X^n$,
with coface maps $d^j:P^{n-1}_{x,y}X \to P^n_{x,y}X$ given by
$$
d^j(x_1,\dots,x_{n-1}) =
\begin{cases}
(x,x_1,\dots,x_{n-1}) & j = 0; \cr
(x_1,\dots,x_j,x_j,\dots,x_{n-1}) & 0<j<n;\cr
(x_1,\dots,x_{n-1},y) & j = n.
\end{cases}
$$

\subsection{Geometric realization} As is well known, each simplicial
topological space $X_\dot$ has a geometric realization $|X_\dot|$, which
is a quotient space
$$
|X_\dot| = \bigg(\coprod_{n\ge 0} X_n \times \D^n\bigg)/\sim
$$
where $\sim$ is a natural equivalence relation generated by identifications
for each morphisms $f : [n] \to [m]$ of $\simp$. If $K$ is an ordered
simplicial complex and $K_\dot$ the associated simplicial set, then $|K_\dot|$
is  homeomorphic to the topological space associated to $K$.

Dually, each cosimplicial space $X\cosmp$ has a kind of geometric realization
$\big\|X\cosmp\big\|$, which is called the {\it total space associated to}
$X^\dot$ (cf.\ \cite{bousfield-kan}). This is exactly the categorical dual
of the geometric realization of a simplicial space. It is simply the subspace
of
$$
\prod_{n\ge 0} X[n]^{\D_n}
$$
consisting of all sequences compatible with all morphisms $f : [n] \to [m]$
in $\simp$, where $X[n]^{\D_n}$ denotes the set of continuous mappings from
$\D_n$ to $X[n]$ endowed with the compact-open topology. Continuous mappings
from a topological space $Z$ to $\big\|X\cosmp\big\|$ correspond naturally to
continuous mappings
$$
\D^\dot \times Z \to X\cosmp
$$
of cosimplical spaces.

As in Section~\ref{it_ints}, we regard $\D^n$ as the time ordered simplex
$$
\D^n = \{(t_1,\dots,t_n) : 0 \le t_1 \le \cdots \le t_n \le 1\}.
$$
There are continuous mappings
$$
PX \to \|P^\dot X\|\text{ and } P_{x,y}X \to \| P_{x,y}^\dot X \|
$$
defined by
$$
\gamma \mapsto
\{(t_1,\dots,t_n) \mapsto (\gamma(0),\gamma(t_1),\dots,\gamma(t_n),\gamma(1))\}
$$
and
$$
\gamma \mapsto \{(t_1,\dots,t_n) \mapsto (\gamma(t_1),\dots,\gamma(t_n))\}
$$
These correspond to the adjoint mappings
$$
\D^\dot \times PX \to P^\dot X \text{ and }
\D^\dot \times P_{x,y}X \to P_{x,y}^\dot X,
$$
which are the continuous mappings of cosimplicial spaces used when defining
iterated integrals in Section~\ref{it_ints}.

\subsection{Cochains}

Applying the singular chain functor to a cosimplicial space $X[\blank]$
yields a simplicial chain complex. Taking alternating sums of the face
maps, we get a double complex $S^\dot(X[\dot];R)$ where
$$
S^{s+t}(X[s];R)
$$
sits in bidegree $(-s,t)$ and total degree $t-s$.\footnote{Note that this
has many elements of negative total degree.} The associated second
quadrant spectral sequence is the Eilenberg-Moore spectral sequence.

Elements of the corresponding total complex can be evaluated on singular
chains $\sigma : \D^t \to \big\| X[\dot]\big\|$ by replacing $\sigma$
by its adjoint
$$
\sigmahat : \D^\dot \times \D^t \to X[\dot].
$$
To evaluate $c \in S^\dot(X[s];R)$ on $\sigma$, first subdivide
$\D^s \times \D^t$ into simplices in the standard way and
then evaluate $c$ on this subdivision of $\D^s\times \D^t \to X[s]$.

When $X$ is a manifold we can apply the de~Rham complex, as above, to obtain
a double complex $E^\dot(X[\dot])$, where $E^s(X[t])$ is placed in bidegree
$(-s,t)$. Integration induces a map of double
complexes
$$
E^\dot(X[\dot]) \to S^\dot(X[\dot];\R).
$$
This is a quasi-isomorphism as is easily seen using the Eilenberg-Moore
spectral sequence.

When $X[\dot]$ is a cosimplical model of a path space, we can say more.
I will treat the case of $P^\dot X$; the case of $P_{x,y}^\dot X$ being
obtained from it by restriction.

The first thing to observe is that $E^\dot(X)^{\otimes(s+2)}$ can be used
in place of
$$
E^\dot(X^{s+2}) = E^\dot(P^s X).
$$
The corresponding double complex has
$$
\big[E^\dot(X)^{\otimes(s+2)}\big]^{s+t}
$$
in bidegree $(-s,t)$. The associated total complex is (essentially by
definition) the {\it unreduced} bar construction
$\B(\E^\dot(X),E^\dot(X),\E^\dot(X))$ on $E^\dot(X)$. Here $E^\dot(X)$ is
considered as a module over itself by multiplication. The chain maps
$$
\B(\E^\dot(X),E^\dot(X),\E^\dot(X)) \to E^\dot(P^\dot X) \to S^\dot(P^\dot X;\R)
$$
are quasi-isomorphisms (use the Eilenberg-Moore spectral sequence). Similarly,
in the case of $P_{x,y}^\dot X$,
\begin{equation}
\label{int_str}
\B(\R,E^\dot(X),\R) \to E^\dot(P_{x,y}^\dot X) \to S^\dot(P_{x,y}^\dot X;\R)
\end{equation}
are quasi-isomorphisms.

We can get cochains on $PX$ by pulling back these along the inclusion
$PX \hookrightarrow \|P^\dot X\|$, which allows us to evaluate elements
of $S^\dot(P^\dot X;R)$ on singular simplices $\sigma : \D^t \to PM$ as
above. In particular, if $\sigma$ is smooth and
$$
w'\otimes w_1 \otimes \dots \otimes w_s \otimes w'' \in
E^\dot(X)^{\otimes(s+2)},
$$
then
\begin{align*}
\langle
\sigma, w' \times w_1 \times w_2 \times \cdots \times w_s \times w''
\rangle
&= \int_{\sigmahat}
\big(w' \times w_1 \times w_2 \times \cdots \times w_s \times w''\big) \cr
&= \pm \big\langle \sigma,
p_0^\ast w' \wedge \bigg(\int w_1 \dots w_r\bigg) \wedge p_1^\ast w''
\big\rangle
\end{align*}
where the sign depends on one's conventions. Thus the cosimplicial
constructions naturally lead to Chen's iterated integrals.

In the case of $P_{x,y}X$, the chain mapping
$B(\R,E^\dot(X),\R) \to \B(\R,E^\dot(X),\R)$
is a quasi-isomorphisms. The cohomology of  $S^\dot(P^\dot_{x,y}X;\Z)$ then
provides the cohomology of iterated integrals with the integral structure
described in Paragraph~\ref{integral} via (\ref{int_str}).

\subsection{Back to Adams}
What is missing from the story so far is chains, which are useful, if not
essential, for computing periods of iterated integrals and mixed Hodge
structures. They are especially useful in situations where the de~Rham theorem
is not true for loop spaces, but where the cohomology of iterated integrals
has geometric meaning. Adams' original work constructs cubical chains on
$P_{x,x} X$ from from certain singular chains on $X$.

Denote the unit interval by $I$ and let $e_j^0, e_j^1 : I^{n-1} \to I^n$
be the $j$th bottom and top face maps of the unit $n$-cube:
$$
e_j^\epsilon :
(t_1,\dots,t_{n-1}) \mapsto (t_1,\dots, t_{j-1},\epsilon,t_j,\dots,t_{n-1}).
$$
For $0\le j \le n$, let $f_j : \D^j \to \D^n$ and $r_j : \D^j \to \D^n$
denote the front and rear $j$-faces of $\D^n$. These correspond to the 
order preserving injections $[j] \to [n]$ uniquely determined by
$f_j(j) = j$ and $r_j(0) = n-j$.

\begin{center}
\begin{figure}[!ht]
\epsfig{file=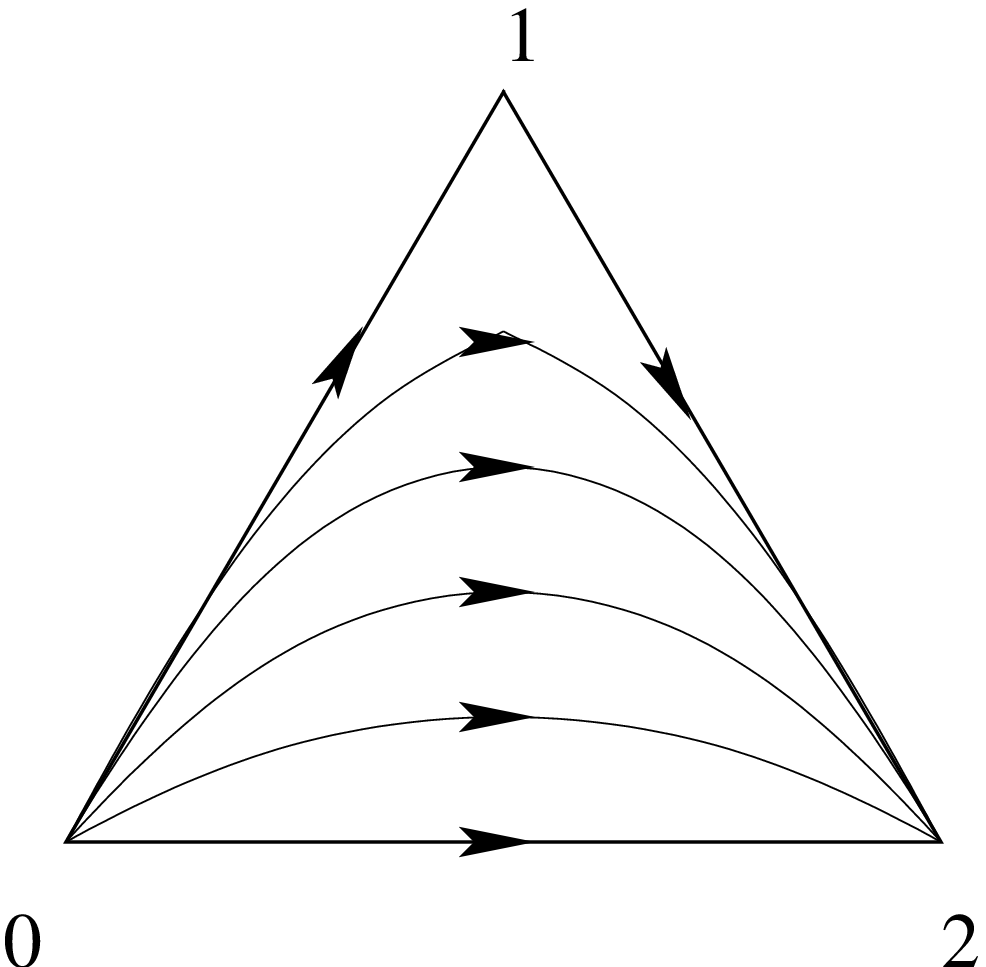, height=1.25in}
\qquad
\epsfig{file=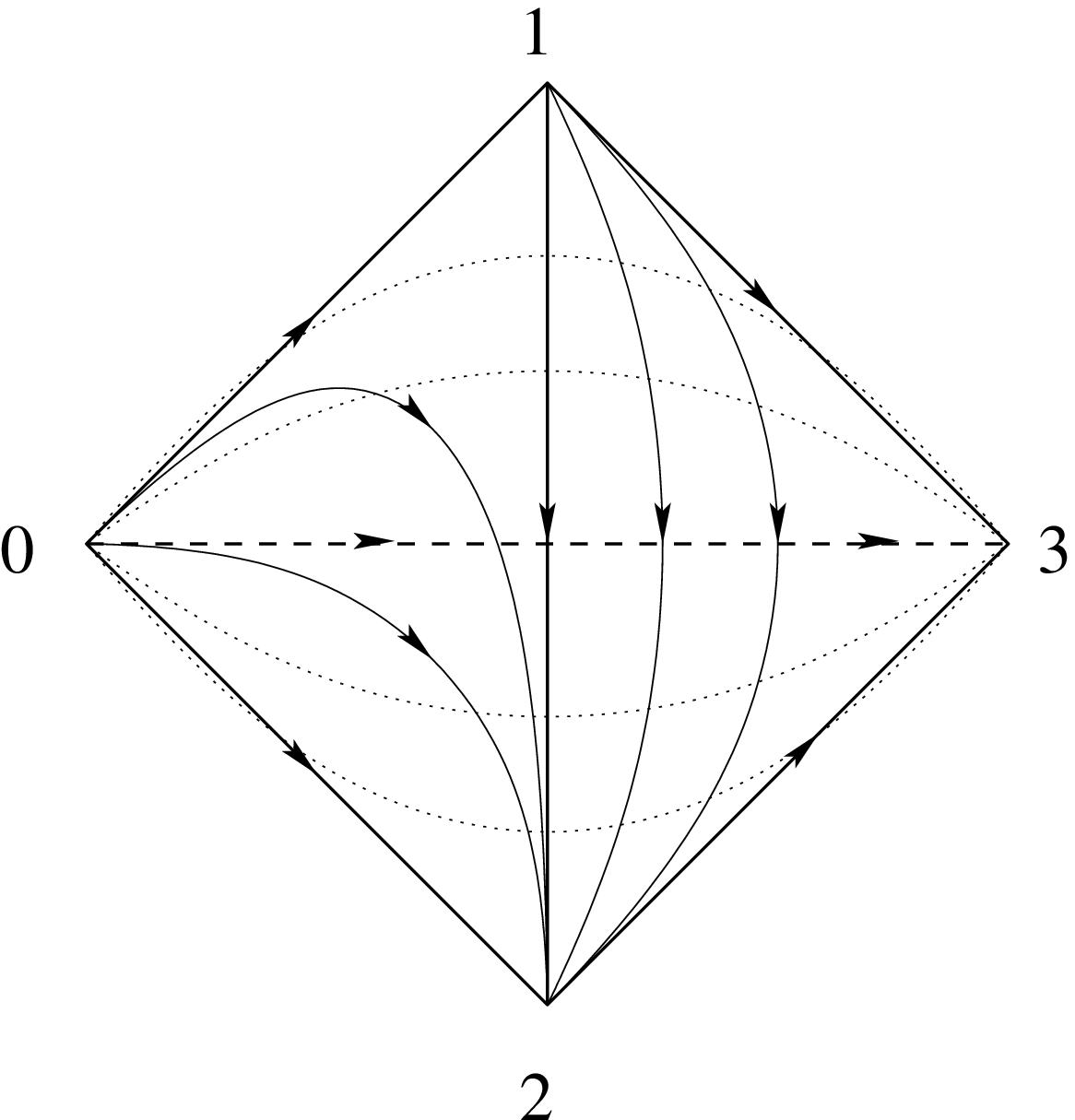, height=1.5in}
\caption{$\theta_2$ and $\theta_3$}\label{simps}
\end{figure}
\end{center}

The starting point is to construct continuous maps\footnote{With care, these
can be made smooth --- details can be found in \cite{chen:bams}.}
$$
\theta_n : I^{n-1} \to P_{0,n} \D^n
$$
with the property that when $0 < j < n$,\footnote{For $\alpha : U \to P_{x,y}X$
and $\beta : V \to P_{y,z} X$, define
$\alpha \ast \beta : U\times V \to P_{x,z}X$ by
$(u,v) \mapsto \alpha(u)\beta(v)$.}
\begin{equation}
\label{cube}
\theta_n \circ e^0_j = P(d^j)\circ\theta_{n-1} : I^{n-2} \text{ and }
\theta_n \circ e^1_j =
\big(P(f_j)\circ\theta_j\big) \ast \big(P(r_{n-j}) \circ \theta_{n-j}\big).
\end{equation}
These are easily constructed by induction on $n$ using the elementary fact
that $P_{0,n}\D^n$ is contractible.
When $n=1$, the unique point of $I^0$ goes  to any path from
$0$ to $1$ in $\D^1$. The cases $n=2$ and 3 are illustrated in
Figure~\ref{simps}.

For a pointed topological space $(X,x)$, let $\So_\dot(X,x)$ be the
subcomplex of the singular chain complex generated by those singular
simplices $\sigma : \D^n \to X$ that map all vertices of $\D^n$ to $x$.
If $X$ is path connected, this computes the integral homology of $X$.

For each such singular simplex $\sigma:\D^n \to X$, we have the singular
cube $P(\sigma)\circ \theta_n : I^{n-1} \to P_{x,x}X$. Set
\begin{equation}
\label{bracket}
[\sigma] =
\begin{cases}
P(\sigma)\circ \theta_n - c_x & n=1; \cr
P(\sigma)\circ \theta_n & n > 1,
\end{cases}
\end{equation}
where $c_x$ denotes the constant loop at $x$. Set\footnote{Strictly
speaking, we need to use Moore paths as we need path multiplication to be
associative.}
$$
[\sigma_1 | \sigma_2 | \dots | \sigma_s] =
[\sigma_1]\ast [\sigma_2]\ast \cdots \ast [\sigma_n]
$$
This extends to an algebra mapping
\begin{align*}
\bigoplus_{s\ge 0} \big(\So_{>0}(X,x)\big)^{\otimes s} &\to
\{\text{reduced cubical chains on }P_{x,x}X \} \cr
\sigma_1 \otimes \cdots \otimes \sigma_s &\mapsto [\sigma_1|\cdots | \sigma_s],
\end{align*}
which is easily seen be to injective. The formula
(\ref{cube}) implies that
\begin{equation}
\label{diff}
\partial [\sigma] = -[\partial \sigma] +
\sum_{1\le j < n} (-1)^j [\sigma_{(j)}|\sigma^{(n-j)}],
\end{equation}
where $\sigma_{(j)}$ denotes the front $j$ face and $\sigma^{(n-j)}$ the
rear $(n-j)$th face of $\sigma$.

Adams' cobar construction is, by definition, the free associative algebra
$$
\Ad(\So_\dot(X,x)) = \bigoplus_{s \ge 0} \big(\So_{>0}(X,x)\big)^{\otimes s}
$$
on $\So_\dot(X,x)$ with the differential (\ref{diff}), where
$\sigma_1 \otimes \dots \otimes \sigma_s$ has degree $-s + \sum \deg \sigma_j$.
This is an augmented, associative, dga, where the augmentation ideal is
generated by the $[\sigma]$.
Adams' main result may be stated by saying that the chain mapping
$$
\Ad(\So_\dot(X,x)) \to \{\text{reduced cubical chains on }P_{x,x}X \}
$$
is a quasi-isomorphism when $X$ is simply connected. Stallings' result
\cite{stallings} for $H_0$ is more elementary.

\begin{proposition}
If $X$ is path connected, then there are natural augmentation preserving
algebra isomorphisms
$$
H_0(\Ad(\So_\dot(X,x))) \cong H_0(P_{x,x}X;\Z) \cong \Z\pi_1(X,x).
$$
\end{proposition}

\begin{proof}[Sketch of proof]
The second isomorphism follows directly from the definitions. We will show
that $H_0(\Ad(\So_\dot(X,x)))$ is isomorphic to $\Z\pi_1(X,x)$. Let
$\Simp_\dot(X,x)$ denote the simplicial set whose $k$-simplices consist of
all singular simplices $\sigma : \D^k \to X$ that map all vertices of $\D^n$
to $x$. After unraveling the definitions (\ref{bracket}) and (\ref{diff}),
we see that $H_0(\Ad(\So_\dot(X,x)))$ is the algebra generated by the
1-simplices $\Simp_1(X,x)$ (augmented by taking the generator corresponding
to each 1-simplex to 1) divided out by the ideal generated by
$\sigma_{01} - \sigma_{02} + \sigma_{12}$,
where $\sigma \in \Simp_2(X,x)$ and $\sigma_{jk}$ is the singular 1-simplex
obtained by restricting $\sigma$ to the edge $jk$ of $\D^2$. It follows from
van Kampen's Theorem that $H_0(\Ad(\So_\dot(X,x)))$ is naturally isomorphic
to the integral group ring of the fundamental group of the geometric
realization of $\Simp_\dot(X,x)$. The result follows as the tautological
mapping $|\Simp_\dot(X,x)| \to X$ is a weak homotopy equivalence.
\end{proof}

With the standard diagonal mapping
$$
\D : \So_\dot(X,x) \to \So_\dot(X,x)\otimes \So_\dot(X,x),
\quad \sigma \mapsto \sum_{0<j<\deg \sigma} \sigma_{(j)} \otimes \sigma^{(n-j)}
$$
$\So_\dot(X,x)$ is a coassociative differential graded coalgebra.
The cobar construction can be defined for any connected, coassociative 
dg coalgebra $C_\dot$. The homology analogue of the  Eilenberg-Moore spectral
sequence implies that if $C_\dot \to \So_\dot(X,x)$ is a dg coalgebra
quasi-isomorphism, then the induced mapping
$$
H_\dot(\Ad(C_\dot)) \to H_\dot(\Ad(\So_\dot(X,x)))
$$
is an isomorphism provided $H_1(X;\Q)=0$, and that
$$
H_0(\Ad(C_\dot))/I^s \to H_0(\Ad(\So_\dot(X,x)))/I^s
$$
is an isomorphism for all $s$ in general.

Elements of $S^\dot(P^\dot_{x,x}X)$ can be evaluated on elements of
$\Ad(\So_\dot(X,x))$ to obtain a chain mapping
$$
S^\dot(P^\dot_{x,x}X;R) \to \Hom(\Ad(\So_\dot(X,x)),R).
$$
which is a quasi-isomorphism for all coefficients $R$ as can be seen
using the Eilenberg-Moore spectral sequence. Consequently, the integral
structures on $H^\dot(\Ch^\dot(P_{x,y}X))$ one obtains from from
$S^\dot(X^\dot;\Z)$ and $\Ad(\So_\dot(X,x))$ agree.

\end{document}